\documentclass[12pt,twoside]{article}
\usepackage{amsmath}
\usepackage{mathrsfs}
\usepackage{amsfonts}
\usepackage{amssymb}

\usepackage{amsfonts,amsmath,amsthm}
\allowdisplaybreaks

\numberwithin{equation}{section} \textwidth 15.5 true cm
\begin{document}
  \title{{\bf Global Low-Energy Weak Solution and Large-Time Behavior for the Compressible Flow of Liquid Crystals
  }}
  \author{  Guochun Wu and Zhong Tan\\{\it\small School of Mathematical Sciences, Xiamen University,
 Fujian 361005, China}\\{\it\small
 Email
 address: guochunwu@126com, ztan85@163.com.}}
  \date{}
  \maketitle
  \begin{abstract}{\small
    We consider the weak solution of the simplified Ericksen-Leslie system modeling compressible nematic liquid crystal flows
    in $\mathbb R^3$. When the initial data is small in $L^2$ and initial density is positive and essentially bounded, we first prove
    the existence of a global weak solution in $\mathbb R^3$. The large-time behavior of a global weak solution is also established.  }
     \\
 \\

  \end{abstract}

\section{Introduction}

$\ \ \ \ $We consider the following hydrodynamic system modeling the
flow of nematic liquid crystal materials [2,7,20],
$$
 \rho_t+\nabla\cdot (\rho u)=0,\eqno(1.1a)$$
$$\rho u_t+\rho u\cdot\nabla u+\nabla P(\rho)=\mu\triangle u+\lambda \nabla div u-\nabla d\cdot\triangle
d,\eqno(1.1b)$$ $$\partial_t d+u\cdot\nabla d=\triangle d +|\nabla
d|^2d,\eqno (1.1c)
$$
for $(t,x)\in[0,+\infty)\times\mathbb{R}^{3}$. Here $\rho$,
$u=(u^{1},u^{2},u^{3})^{t}$ and $P$ denote the density, the
velocity, and the pressure respectively. $d=(d^1,d^2,d^3)^t$ is the
unit-vector ($|d|=1$) on the sphere $\mathbb S^2\subset \mathbb R^3$
representing the macroscopic molecular orientation of the liquid
crystal materials. $\mu$ and $\lambda$ are positive viscosity
constants, and $div$ and $\triangle$ are the usual spatial
divergence and Laplace operators.

The above system (1.1) is a simplified version of the
Ericksen-Leslie model for the hydrodynamics of nematic liquid
crystals. The mathematical analysis of the incompressible liquid
crystal flows was initialed by Lin and Liu in [21,22]. In [29,30],
Wang proved the global existence of strong solutions in whole space
under some small conditions. In dimension two, the global existence
and uniqueness of weak solutions were studied in [10,11,18,23,31]
and reference therein. On the other hand, the global existence of
weak solutions to incompressible liquid crystal in $\mathbb R^3$ is
still an outstanding open question.

When the fluid is allowed to be compressible, the Ericksen-Leslie
system becomes more complicate. To our knowledge, there seems very
few analytic works available yet. The local-in-time strong solutions
to the initial value or initial boundary  value problem of system
(1.1) with nonnegative initial density were studied in [3,13]. Based
on [16,17], the blow up criterion of strong solutions were obtained
in [13,14]. The global existence and uniqueness of strong solution
in critical space were studied in [12]. Motivated by [15], when the
initial data were sufficiently smooth and were suitably small in
some energy-norm, the global well-posedness of classical solutions
were proved in [19]. Especially global weak solutions  in two
dimensions was established in [6] under some small condition or
geometric angle condition.

Our aim in this paper is to establish the global existence of
low-energy weak solutions of system (1.1), with the following
initial conditions:
$$(\rho(\cdot,0),u(\cdot,0),d(\cdot,0))=(\rho_0,u_0,d_0),\eqno(1.2)$$
where $\rho_0$ is bounded above and below away from zero, $|d_0|=1$,
$u_0,\nabla d_0\in L^p(\mathbb R^3)$ for some $p>6$, and modulo
constants, $(\rho_0,u_0,d_0)$ is small in $L^2(\mathbb R^3)$. Thus
the total initial energy is small, but no other smallness or
regularity conditions are imposed.

When the direction field $d$ does not appear, (1.1) reduces to the
compressible Navier-Stokes equations. The global classical solutions
were first obtained by Matsumura-Nishida [25,26] for initial data
close to a non-vacuum equilibrium in $H^3(\mathbb R^3)$. In
particular, the theory requires that the solution has small
oscillations from a uniform non-vacuum state so that the density is
strictly away from the vacuum and the gradient of the density
remains bounded uniformly in time. Later , Hoff [8,9] studied the
problem for discontinuous initial data. For the existence of
solutions for arbitrary data, the major breakthrough is due to Lions
[24] (see also Feireisl [5]), where he obtains global existence of
weak solutions-defined as solutions with finite energy. Suen and
Hoff [28] adopted Hoff's techniques to obtain global existence of
low-energy weak solutions for the magnetohydrodynamics. In this
paper, we shall study the Cauchy problem (1.1)-(1.2) for liquid
crystals and establish the global existence and large time behavior
of low-energy weak solutions However, compared with the compressible
Navier-Stokes equations, some new difficulties arise due to the
additional presence of the liquid crystal directional field.
Especially, the super critical nonlinearity $|\nabla d|^2d$ in the
transported heat flow of harmonic map equation $(1.1c)$ and the
strong coupling nonlinearity $\nabla d\cdot\triangle d$ in the
momentum equations $(1.1b)$ will cause serious difficulties in the
proofs of the time-independent global energy estimates.

To state the main results in a precise way, we first introduce some
notations and conventions which will be used throughout the paper.
For a given unit vector $n\in\mathbb S^2$ and a positive integer m,
we denote
$$H_n^m(\mathbb R^3;\mathbb S^2):=\{d:d-n\in H^m(\mathbb R^3),|d|=1\
a.e.\ in\ \mathbb R^3\}.$$ We use the usual notation for H\" older
seminorms: for $v:\mathbb R^3\longrightarrow \mathbb R^m$ and
$\alpha\in (0,1]$,
$$\langle v\rangle^\alpha=\sup\limits_{x_1, x_2\in\mathbb R^3\atop x_1\neq x_2}\frac{|v(x_2)-v(x_1)|}{|x_2-x_1|^\alpha};$$
and for $v:Q\subseteq\mathbb
R^3\times[0,\infty)\longrightarrow\mathbb R^3$ and
$\alpha_1,\alpha_2\in(0,1]$,
$$\langle v\rangle^{\alpha_1,\alpha_2}_{Q}=\sup\limits_{(x_1,t_1),( x_2,t_2)\in Q
\atop (x_1,t_1)\neq(
x_2,t_2)}\frac{|v(x_2,t_2)-v(x_1,t_1)|}{|x_2-x_1|^{\alpha_1}+|t_2-t_1|^{\alpha_2}}.$$
If $X$ is a Banach space we will abbreviate $X^3$ by $X$ when
convenient. Finally if $I\subset[0,\infty)$ is an interval,
$C^1(I;X)$ will be the elements $v\in C(I;X)$ such that the
distribution derivative $v_t\in \mathcal {D}'(int\ I;\mathbb R^3)$
is realized as an element of $C(I;X)$.

As it was pointed out in [8], the effective viscous flux plays an
important role in the mathematical theory of compressible fluid
dynamics. More precisely, let $F$ and $\omega$ be the effective flux
and vorticity defined by
$$F\triangleq (\mu+\lambda)divu-(P(\rho)-P(\tilde{\rho}))\ \ and \ \ \omega\triangleq\nabla\times u.\eqno(1.3)$$
It is not hard to check that
$$\nabla d\triangle d=div(\nabla d\otimes \nabla d)-\nabla\frac{|\nabla d|^2}{2}.$$
So, it follows from (1.1b) that
$$\triangle F=div(\rho\dot{u}+div(\nabla d\otimes \nabla d)-\nabla\frac{|\nabla d|^2}{2}),\ \ \ \mu\triangle\omega
=\nabla\times(\rho\dot u+div(\nabla d\otimes \nabla d)),\eqno(1.4)$$
where $``\ \dot\ \ "$ denotes the material derivative, i.e.,
$$\dot f:=\partial_tf+u\cdot\nabla f.$$

Now we give a precise formulation of our results. First concerning
the pressure P, we focus our interest on the case of isentropic
flows and assume that
$$P(\rho)=a\rho^\gamma\ \ with\ \ a>0,\gamma\geq 1.\eqno(1.5)$$
Next we fix a positive reference density $\tilde{\rho}$ and then
choose positive bounding densities $\underline{\rho}$ and
$\bar{\rho}$ satisfying
$$\underline \rho<\min\{\tilde \rho,\rho'\}\ \ and \ \ \max\{\tilde\rho,\rho''\}<\bar \rho,\eqno(1.6)$$
and finally we define a positive number $\delta$ by
$$\delta=\min\{\min\{\tilde\rho,\rho'\}-\underline\rho,\bar\rho-\max\{\tilde\rho,\rho''\},\frac{1}{2}(\bar\rho-\underline\rho)\}.\eqno(1.7)$$
(Notice that $\delta$ need not be ``small" in the usual sense.)
Concerning the diffusion coefficients $\mu$ and $\lambda$ we assume
that
$$0\leq\lambda<\frac{3+\sqrt{21}}{6}\mu.\eqno(1.8)$$
It follows that
$$\frac{1}{4}\mu(p-2)-\frac{[\frac{1}{4}\lambda(p-2)]^2}{\frac{1}{3}\mu+\lambda}>0\eqno(1.9)$$
for $p=6$ and consequently for some $p>6$, which we now fix.

Concerning the initial data $(\rho_0,u_0,d_0)$, we assume there is a
positive number $N$, which may be arbitrary large, and a positive
number $b<\delta$ such that
$$\|u_0\|_{L^p}+\|\nabla d_0\|_{L^p}\leq N\eqno(1.10)$$
and
$$\underline \rho+b<ess\inf\rho_0\leq ess\sup\rho_0<\bar\rho-b.\eqno(1.11)$$
We assume also that $$d_0-n\in H_n^1(\mathbb R^3;\mathbb
S^2)\eqno(1.12)$$ and write
$$C_0\triangleq \int_{\mathbb R^3}(\frac{1}{2}\rho_0|u_0|^2+G(\rho_0)+\frac{1}{2}|\nabla d_0|^2)dx,\eqno(1.13)$$
where $G(\rho)$ is the potential energy density defined by
$$G(\rho)=\rho\int^\rho_{\tilde{\rho}}\frac{P(s)-P(\tilde\rho)}{s^2}ds.\eqno(1.14)$$
It is clear that there exist two positive constant $c_1,c_2$ only
depending on $\underline\rho,\bar\rho,$ and $\tilde\rho$
$$c_1(\underline\rho,\bar\rho,\tilde\rho)(\rho-\tilde\rho)^2\leq G(\rho)\leq c_2(\underline\rho,\bar\rho,\tilde\rho)(\rho-\tilde\rho)^2.\eqno(1.15)$$

Weak solutions of (1.1)-(1.2) are defined in a usual way.\\
\\{\bf Definition 1.1.} A pair of functions $(\rho,u,H)$ is said to
be a weak solution of (1.1)-(1.2) provided that
$(\rho-\tilde\rho,\rho u)\in C([0,\infty);H^{-1}(\mathbb R^3))$,
$d-n\in C([0,\infty);L^{2}(\mathbb R^3))$, $u,\nabla d\in
L^\infty([0,\infty);L^2(\mathbb R^3))$, $\nabla u\in
L^2((0,\infty);L^2(\mathbb R^3))$, and $|d(\cdot,t)|=1\ a.e.\ in\
\mathbb R^3$ for $t\geq 0$. Moreover, the following identities hold
for any test function $\psi\in\mathcal {D}(\mathbb
R^3\times(t_1,t_2))$ with $t_2\geq t_1\geq 0$ and $j=1,2,3$:
$$\int_{\mathbb R^3}\rho\psi(x,t)dx|^{t_2}_{t_1}=\int^{t_2}_{t_1}\int_{\mathbb R^3}(\rho\psi_t+\rho u\cdot\nabla\psi)dxdt,$$
$$\begin {array}{rl}&\int_{\mathbb R^3}\rho u^j\psi(x,t)dx|^{t_2}_{t_1}+\int^{t_2}_{t_1}\int_{\mathbb R^3}(\mu\nabla u^j\cdot\nabla\psi+\lambda(div u)\psi_{x_j})dxdt\\=&
\int^{t_2}_{t_1}\int_{\mathbb R^3}(\rho u^j\psi_t+\rho
u^ju\cdot\nabla\psi+P(\rho)\psi_{x_j}-\frac{1}{2}|\nabla
d|^2\psi_j+d_{x_j}\nabla d\cdot\nabla\psi)dxdt,\end {array}$$ and
$$\begin {array}{rl}&\int_{\mathbb R^3}(d^j-n^j)(x,t)\psi(x,t)dx|^{t_2}_{t_1}+\int^{t_2}_{t_1}\int_{\mathbb R^3}\nabla d^j\nabla\psi dxdt\\=&
\int^{t_2}_{t_1}\int_{\mathbb R^3}((d^j-n^j)\psi_t-u\cdot\nabla
d^j\psi+|\nabla d|^2d^j\psi)dxdt.\end {array}$$

Our main results are formulated as the following theorem.\\
\\{\bf Theorem 1.2.} Assume that the system parameters in (1.1)
satisfy the conditions (1.5)-(1.9) and let positive numbers $N$ and
$b<\delta$ be given. Then there are positive constants
$\varepsilon$, $C$, and $\theta$ depending on the parameters and
assumptions in (1.5)-(1.9), on $N$, and on a positive lower bound
for $b$, such that, if initial data $(\rho_0,u_0,d_0)$ are given
satisfying (1.10)-(1.14) with
$$C_0<\varepsilon,\eqno(1.16)$$
then there is a solution $(\rho,u,d)$ to (1.1)-(1.2) in the sense of
the definition 1.1. Moreover, the solution satisfies the following:
$$\rho-\tilde\rho,\rho u\in C([0,\infty);H^{-1}(\mathbb R^3)),\eqno(1.17)$$
$$d-n\in C([0,\infty);L^{2}(\mathbb R^3)),\eqno(1.18)$$
$$\nabla u\in L^2([0,\infty);\mathbb R^3),\eqno(1.19)$$
$$u(\cdot,t),\nabla d(\cdot,t)\in H^1(\mathbb R^3),\ \ t>0\eqno(1.20)$$
$$F(\cdot,t),\omega(\cdot,t)\in H^1(\mathbb R^3),\ \ t>0\eqno(1.21)$$
$$\langle u\rangle^{\frac{1}{2},\frac{1}{8}}_{\mathbb R^3\times[\tau,\infty)},
\langle \nabla d\rangle^{\frac{1}{2},\frac{1}{8}}_{\mathbb
R^3\times[\tau,\infty)}\leq C(\tau)C_0^\theta,\eqno(1.22)$$ where
$C(\tau)$ may depend additionally on a positive lower bound for
$\tau$,
$$\underline \rho\leq\rho(x,t)\leq\bar\rho\ a.e.\ on\ \mathbb R^3\times[0,\infty),\eqno(1.23)$$
and
$$\begin {array}{rl} &\sup\limits_{t>0}\int_{\mathbb R^3}[|\rho-\tilde\rho|^2+|u|^2+\frac{1}{2}|\nabla d|^2+\sigma(|\nabla u|^2+|\nabla^2d|^2)
+\sigma^5(F^2+|\nabla\omega|^2)]dx\\&+\int_0^\infty\int_{\mathbb
{R}^3}[|\nabla u|^2+|\triangle d+|\nabla d|^2d|^2+\sigma(|\dot
u|^2+|\nabla d_t|^2+|\nabla \omega|^2)\\&+\sigma^5(|\nabla\dot
u|^2+|\nabla^2 d_t|^2)]dxds\\&\leq CC_0^\theta,
\end {array}\eqno(1.24)$$
where $\sigma(t)=\min\{1,t\}$. Moreover we also have the following
large-time behavior:
$$\lim_{t\rightarrow\infty}(\|\rho-\tilde{\rho}\|_{L^l(\mathbb R^3)}+\|u\|_{W^{1,r}(\mathbb R^3)}+\|\nabla d\|_{W^{1,r}(\mathbb R^3)})=0,\eqno(1.25)$$
holds for $l\in(2,\infty),r\in (2,6)$.\\

The rest of this paper is devoted to prove Theorem 1.2. In Section
2, we collect some useful inequalities and basic results. In Section
3, we derive the time-independent energy estimates of the solution.
The key pointwise upper and lower bound of the density are
established in Section 4. Finally, the proof of the main results
will be done in Section 5.

\section{Preliminaries}

$\ \ \ \ $In this section, we state some auxiliary lemmas, which
will be frequently used in the sequel. We start with the well-known
Gagliardo-Nirenberg inequality (see, for instance, [1,32]).\\
\\{\bf Lemma 2.1.} First, given $r\in [2,6]$ there is a constant
$C(r)$ such that for $f\in H^1(\mathbb R^3)$,
$$\|f\|_{L^r(\mathbb R^3)}\leq C(r)\|f\|^{(6-r)/2r}_{L^2(\mathbb R^3)}\|\nabla f\|^{(3r-6)/2r}_{L^2(\mathbb R^3)}.\eqno(2.1)$$
Next, for any $r\in (3,\infty)$ and $q>1$, there is a constant
$C(r,q)$ such that for $f\in L^q(\mathbb R^3)\cap W^{1,r}(\mathbb
R^3)$,
$$\|f\|_{L^\infty(\mathbb R^3)}\leq C(r,q)\| f\|^{q(r-3)/(3r+q(r-3))}_{L^q(\mathbb R^3)}\|\nabla f\|^{3r/(3r+q(r-3))}_{L^r(\mathbb R^3)},\eqno(2.2)$$
and
$$\langle f\rangle^\alpha_{\mathbb R^3}\leq C(r)\|\nabla f\|_{L^{r}(\mathbb R^3)},\eqno(2.3)$$
where $\alpha=1-\frac{3}{r}$.\\

The next lemma is due to Hoff [8], which will be used to prove the
uniform (in time) bound of density.\\
\\{\bf Lemma 2.2.} If $\Gamma$ is the fundamental solution for the
Laplace operator in $\mathbb R^3$, then given $p_1\in [1,3)$ and
$p_2\in(3,\infty]$, there is a constant $C=C(n,p_1,p_2)$ such that
$$\|\Gamma_{x_j}*f\|_{L^\infty(\mathbb R^3)}\leq C(n,p_1,p_2)[\|f\|_{L^{p_1}(\mathbb R^3)}+\|f\|_{L^{p_2}(\mathbb
R^3)}].\eqno(2.4)$$\\

Finally, we need the local-in-time existence theorem of (1.1)-(1.2).
Using the mollifier technique, the local solutions which can be
proved rigorously by the standard method of Matsumura and Nishida
[25,26].\\
\\{\bf Proposition 2.3.} Assume that the initial data $(\rho_0,u_0,d_0)$
satisfies
$$(\rho_0-\tilde\rho,u_0)\in H^3(\mathbb R^3)\ d_0\in H_n^4(\mathbb R^3;\mathbb S^2)\ and\ \inf_{x\in\mathbb R^3}(\rho_0(x))>0.\eqno(2.5)$$
Then there exists a positive time $T_0$, which may depend on
$\inf_{x\in\mathbb R^3}(\rho_0(x))$, such that the Cauchy problem
(1.1)-(1.2) has a unique smooth solution $(\rho,u,d)$ on $\mathbb
R^3\times[0,T_0]$ satisfying
$$\rho(x,t)>0\ for\ all\ x\in\mathbb R^3, t\in [0,T_0],\eqno(2.6)$$
$$\rho-\tilde\rho\in C([0,T_0];H^3)\cap C^1([0,T_0];H^2),\eqno (2.7)$$
$$u\in C([0,T_0];H^3)\cap C^1([0,T_0];H^1),\eqno (2.8)$$
$$d-n\in C([0,T_0];H_n^4(\mathbb R^3;\mathbb S^2))\cap C^1([0,T_0];H_n^2(\mathbb R^3;\mathbb S^2)).\eqno
(2.9)$$

In view of Lemma 2.1 and the classical estimates of elliptic system, we have\\
\\{\bf Lemma 2.4.} Let $(\rho,u,d)$ be as in Proposition 2.3, then there exists a generic positive constant $C$,
depending only on $\mu$ and $\lambda$ and $p$, such that for
$p\in(1,\infty)$
$$\|\nabla F\|_{L^p(\mathbb R^3)}+\|\nabla\omega\|_{L^p(\mathbb R^3)}\leq C(\|\rho\dot u\|_{L^p(\mathbb R^3)}+\|\nabla d\cdot\triangle d\|_{L^p(\mathbb R^3)}),\eqno(2.10)$$
$$\|\nabla u\|_{L^p(\mathbb R^3)}\leq C(\|F\|_{L^p(\mathbb R^3)}+\|\omega\|_{L^p(\mathbb R^3)}+\|P(\rho)-P(\tilde\rho)\|_{L^p(\mathbb R^3)}),\eqno(2.11)$$
where $F$ and $\omega$ are defined in (1.4).\\
\\{\bf Proof.} An application of the $L^p$-estimate of elliptic
systems to (1.4) gives (2.10). On the other hand, since $-\triangle
u=-\nabla div u+\nabla\times \omega$, it holds that
$$\nabla u=-\nabla(-\triangle)^{-1}\nabla div u+\nabla(-\triangle)^{-1}\nabla\times\omega,$$
which, combined with the Marcinkiewicz multiplier theorem (Stein
[20,p.96]), we arrive at
$$\begin {array} {rl}\|\nabla u\|_{L^p(\mathbb R^3)}&\leq C(\|div u\|_{L^p(\mathbb R^3)})+\|\omega\|_{L^p(\mathbb R^3)}
\\&\leq C(\|F\|_{L^p(\mathbb R^3)}+\|\omega\|_{L^p(\mathbb R^3)}+\|P(\rho)-P(\tilde\rho)\|_{L^p(\mathbb R^3)}).\end {array}$$
Thus the proof the lemma is completed.

\section{A priori estimates}

$\ \ \ \ $This section is devoted to establish a number of a priori
bounds for local-in-time smooth solutions, corresponding roughly to
(1.24). Those are rather long and technical. We have therefore
omitted those which are identical to or nearly identical to
arguments given elsewhere in the literature of whose details we
regard as routine. On the other hand, we have endeavored to describe
the flow of the arguments in such a way that the diligent reader can
reconstruct the details without undue difficulty.

Let $T>0$ be fixed and assume that $(\rho,u,d)$ is a smooth solution
of (1.1)-(1.2). We define a functional $A(t)$ for a given such
solution that
$$\begin {array}{rl}A(t)=&\sup\limits_{0\leq s\leq t }\int_{\mathbb R^3}[\sigma(|\nabla u|^2+|\nabla^2 d|^2)+\sigma^5(|\dot u|^2+|\nabla\omega|^2+|\nabla d_t|^2)]dx
\\&+\int_0^t\int_{\mathbb
{R}^3}[\sigma(|\dot u|^2+|\nabla d_t|^2+|\nabla
\omega|^2)+\sigma^5(|\nabla\dot u|^2+|\nabla^2 d_t|^2)]dxds,\end
{array}\eqno(3.1)$$ where $\sigma(t)=\min\{1,t\}$, and we obtain the
following a priori bound for $A(t)$ under the assumptions that the
initial energy $C_0$ in (1.13) is small enough and that the density
remains bounded above and below away from zero:\\
\\{\bf Proposition 3.1.} Assume that the system parameters in (1.1)
satisfy the conditions in (1.5)-(1.8) and let positive numbers $N$
and $b<\delta$ be given. Assume $(\rho,u,d)$ is a solution of (1.1)
on $\mathbb R^3\times [0,T]$ in the sense of Proposition 2.3 with
initial data $(\rho_0,u_0)\in H^3(\mathbb R^3)$ and $d_0\in
H_n^4(\mathbb R^3;\mathbb S^2)$ satisfying (1.10)-(1.13), then there
are positive constants $\varepsilon,M,$ and $\theta$ depending on
the parameters and assumptions in (1.5)-(1.8), on $N$, and a
positive lower bound for $b$, such that if  $C_0<\varepsilon$ and
$$\underline\rho\leq\rho(x,t)\leq \bar\rho\ on\ \mathbb R^3\times[0,T],$$
then
$$A(T)\leq MC_0^\theta.$$

The proof will be given in a sequence of lemmas in which we estimate
a number of auxiliary functionals. To describe these we first recall
the definition of (1.9) of $p$, which is an open condition, and
which therefore allows us to choose $q\in[6,\min\{p,12\})$ which
also satisfies (1.9). Then for given $(\rho,u,d)$ we define
$$\begin {array}{rl}\bar A(t)=&\sup\limits_{1\leq s\leq t }\int_{\mathbb R^3}(|\nabla u|^2+|\nabla^2 d|^2+|\dot u|^2+|\nabla\omega|^2+|\nabla d_t|^2)dx
\\&+\int_0^t\int_{\mathbb
{R}^3}(|\dot u|^2+|\nabla d_t|^2+|\nabla \omega|^2+|\nabla\dot
u|^2+|\nabla^2 d_t|^2)dxds,\end {array}$$
$$\begin {array}{rl} B_q(t)=&\sup\limits_{0\leq s\leq t }\int_{\mathbb R^3}(| u|^q+|\nabla d|^q)dx+\int_0^t\int_{\mathbb
{R}^3}(| u|^{q-2}|\nabla u|^2+|\nabla d|^{q-2}|\nabla^2d|^2)dxds
\\&+\int_0^t\int_{\mathbb {R}^3}(|
u|^{q-4}|\nabla (|u|^2)|^2+|\nabla d|^{q-4}|\nabla (|\nabla
d|^2)|^2)dxds,\end {array}$$
$$D(t)=\sup_{0<s\leq t}\int_{\mathbb R^3}\sigma^5(|\nabla d|^2|\nabla^2 d|^2+| u|^2|\nabla^2 d|^2+|\nabla u|^2|\nabla d|^2)dx,$$
$$\bar D(t)=\sup_{1\leq s\leq t}\int_{\mathbb R^3}(|\nabla d|^2|\nabla^2 d|^2+| u|^2|\nabla^2 d|^2+|\nabla u|^2|\nabla d|^2)dx,$$
$$\begin {array} {rl}E(t)=&\int_0^{ t}\int_{\mathbb R^3}[\sigma^\frac{3}{2}(|\nabla u|^3+|\nabla^2d|^3)+\sigma^5(|\nabla u|^4+|\nabla ^2d|^4)]dxds
\\&+|\sum\limits_{1\leq k_i,j_m\leq 3}\int_0^t\int_{\mathbb R^3}\sigma u^{j_1}_{x_{k_1}}u^{j_2}_{x_{k_2}}u^{j_3}_{x_{k_3}}dxds|,\end {array}$$
and
$$\bar E(t)=\int_1^t\int_{\mathbb R^3}(|\nabla u|^3+|\nabla d|^2|\nabla^2d|^2+|\nabla u|^4+|\nabla ^2d|^4)dxds.$$

It will be seen that the assumed regularity (2.6)-(2.9) suffices to
justify the estimates that follow. We begin with the following $L^2$
energy estimate.\\
\\{\bf Lemma 3.2.} Assume that the hypotheses and notations of
Proposition 3.1 are in force. Then
$$\begin {array}{rl} &\sup\limits_{0\leq t\leq T}\int_{\mathbb R^3}(|\rho-\tilde{\rho}|^2+|u|^2+|\nabla d|^2)dx\\&\ \ \ +\int_0^T\int_{\mathbb
{R}^3}(|\nabla u|^2+|\triangle d+|\nabla d|^2d|^2)dxdt\leq MC_0.
\end {array}\eqno(3.2)$$\\
\\{\bf Proof.} Multiplying (1.1b) by $u$ and integrating over $\mathbb
R^3$, we have
$$\frac{d}{dt}\int_{\mathbb R^3}\frac{1}{2}\rho|u^2|dx+\int_{\mathbb R^3}\nabla P(\rho)udx+\int_{\mathbb R^3}(\mu|\nabla
u|^2+\lambda|divu|^2)dx=- \int_{\mathbb R^3}u\cdot\nabla
d\cdot\triangle d.\eqno(3.3)$$ By the mass equation (1.1a) and the
definition of $G(\rho)$ in (1.14), we have
$$G(\rho)_t+div(G(\rho)u)+(P(\rho)-\tilde P)div\ u=0.$$
Integrating and adding the result to (3.3) we obtain
$$\frac{d}{dt}\int_{\mathbb R^3}\frac{1}{2}\rho|u^2|+G(\rho)dx+\int_{\mathbb R^3}(\mu|\nabla
u|^2+\lambda|divu|^2)dx=- \int_{\mathbb R^3}u\cdot\nabla
d\cdot\triangle d.\eqno(3.4)$$ Multiplying (1.1c) by $\triangle
d+|\nabla d|^2d$ and integrating over $\mathbb R^3$, using
integration by parts and the fact that $|d|=1$ we obtain
$$\frac{d}{dt}\int_{\mathbb R^3}\frac{1}{2}|\nabla d|^2dx+\int_{\mathbb
{R}^3}|\triangle d+|\nabla d|^2d|^2dx=\int_{\mathbb R^3}u\cdot\nabla
d\cdot\triangle d.\eqno(3.5)$$ Adding (3.4) to (3.5) and integrating
over [0,t], yields (3.2) by (1.15). Thus the proof of lemma is completed.\\
\\{\bf Lemma 3.3.} Assume that the hypotheses and notations of
Proposition 3.1 are in force. Then for $0<t\leq 1\wedge T$
$$\begin {array}{rl}&\sup\limits_{0<s\leq t}\sigma\int_{\mathbb R^3}(|\nabla u|^2+|\nabla^2d|^2)dx+\int_0^t\int_{\mathbb R^3}
\sigma(|\dot u|^2+|\nabla d_t|^2)dxds\\ \leq&
M[C_0+C_0^{\frac{q-4}{q-2}}B_q^{\frac{2}{q-2}}+C_0^{\frac{q-6}{q-2}}B_q^{\frac{4}{q-2}}
+E],\end {array}\eqno (3.6)$$ and if $T>1$ and $1\leq t\leq T$, then
$$\begin {array}{rl}&\sup\limits_{1\leq s\leq t}\int_{\mathbb R^3}(|\nabla u|^2+|\nabla^2d|^2)dx+\int_0^t\int_{\mathbb R^3}
(|\dot u|^2+|\nabla d_t|^2)dxds\\ \leq& M[C_0+C_0^\frac{3}{2}\bar
A^\frac{1}{2})+\bar
E]+A(1).\end {array}\eqno (3.7)$$ \\
{\bf Proof.} For $0\leq t\leq 1\wedge T$, multiplying the equation
(1.1b) by $\sigma\dot u$ and integrating over $\mathbb R^3\times
[0,t]$, we have
$$\begin {array}{rl}&\sup\limits_{0<s\leq t}\sigma\int_{\mathbb R^3}|\nabla u|^2dx+\int_0^t\int_{\mathbb R^3}
|\dot u|^2dxds\\ \leq& M\{C_0+|\int_0^t\int_{\mathbb R^3}\sigma[\dot
u(div(\nabla d\otimes\nabla d)-\frac{1}{2}\nabla|\nabla
d|^2)dxds|+E\}.\end {array}\eqno(3.8)$$ Differentiating (1.1c) with
respect to x, we have
$$\nabla d_t-\triangle\nabla d=\nabla(|\nabla d|^2d-u\cdot\nabla d).\eqno(3.9)$$
Multiplying the above equation by $\sigma\nabla d_t$ and integrating
over $\mathbb R^3\times [0,t]$, we have
$$\begin {array}{rl}&\frac{1}{2}\sigma\int_{\mathbb R^3}|\triangle d|^2dx+\int_0^t\int_{\mathbb R^3}\sigma|\nabla d_t|^2dxds\\
=&\frac{1}{2}\int_0^t\int_{\mathbb R^3}\sigma'| \triangle
d|^2dxds+\int_0^t\int_{\mathbb R^3}\sigma\nabla d_t\nabla(|\nabla
d|^2d-u\cdot\nabla d)dxds.
\end {array}$$
Adding this to (3.8) and combining with Cauchy's inequality we then
get
$$\begin {array}{rl}&\sup\limits_{0<s\leq t}\sigma\int_{\mathbb R^3}(|\nabla u|^2+|\triangle d|^2)dx+\int_0^t\int_{\mathbb R^3}
\sigma(|\dot u|^2+|\nabla d_t|^2)dxds\\ \leq&
M\{C_0+E+\int_0^t\int_{\mathbb R^3}\sigma'| \triangle
d|^2dxds\\&+\int_0^t\int_{\mathbb R^3}[\sigma(|\nabla^2 d|^2|\nabla
d|^2+|\nabla u|^2|\nabla d|^2+|\nabla^2 d|^2|u|^2)]dxds\}.\end
{array}\eqno(3.10)$$ The right terms can be estimated as follows:
$$\begin {array}{rl}\int_0^t\int_{\mathbb R^3}\sigma'| \triangle
d|^2dxds&\leq M\int_0^{t\wedge 1}\int_{\mathbb R^3}| \triangle
d+|\nabla d|^2d|^2+|\nabla d|^4dxds\\&\leq M[C_0+(\int_0^{t\wedge
1}\int_{\mathbb R^3}|\nabla
d|^2dxds)^{\frac{q-4}{q-2}}\int_0^{t\wedge 1}\int_{\mathbb
R^3}|\nabla d|^qdxds)^{\frac{2}{q-2}}]\\&\leq
M(C_0+C_0^{\frac{q-4}{q-2}}B_q^{\frac{2}{q-2}}),\end {array}$$
$$\begin {array}{rl}\int_0^t\int_{\mathbb R^3}\sigma|\nabla u|^2|\nabla d|^2dxds&\leq M
(\int_0^t\int_{\mathbb R^3}\sigma^{\frac{3}{2}}|\nabla
u|^3dxds)^\frac{2}{3}(\int_0^t\int_{\mathbb R^3}|\nabla
d|^6dxds)^\frac{1}{3}\\&\leq M[H+(\int_0^t\int_{\mathbb R^3}|\nabla
d|^2dxds)^\frac{q-6}{q-2}(\int_0^t\int_{\mathbb R^3}|\nabla
d|^qdxds)^\frac{4}{q-2}]\\&\leq
M(C_0+C_0^{\frac{q-6}{q-2}}B_q^{\frac{4}{q-2}}).
\end {array}$$
The other two terms in the integral on the right side of (3.10) are
bounded in a similar way, and (3.6) follows.

For $1\leq t\leq T$, as in (3.10), we have
$$\begin {array}{rl}&\sup\limits_{1\leq s\leq t}\int_{\mathbb R^3}(|\nabla u|^2+|\triangle d|^2)dx+\int_1^t\int_{\mathbb R^3}
(|\dot u|^2+|\nabla d_t|^2)dxds\\ \leq& M\{C_0+\bar
E+\int_1^t\int_{\mathbb R^3}[(|\nabla^2 d|^2|\nabla d|^2+|\nabla
u|^2|\nabla d|^2+|\nabla^2 d|^2|u|^2)]dxds\}+A(1).\end
{array}\eqno(3.11)$$ Using the fact that
$$|\nabla d|^2=-d\cdot\triangle d\ (since \ |d|=1),\eqno(3.12)$$
the right terms can be bounded as follows:
$$\begin {array}{rl}\int_1^t\int_{\mathbb R^3}|\nabla u|^2|\nabla d|^2dxds&\leq M\int_1^t\int_{\mathbb R^3}|\nabla u|^3+|\nabla d|^6dxds\\&\leq M
\int_1^t\int_{\mathbb R^3}|\nabla u|^3+|\nabla d|^2|\nabla^2
d|^2dxds\\&\leq M\bar E\end {array}$$
$$\begin {array}{rl}\int_1^t\int_{\mathbb R^3}|\nabla^2 d|^2|u|^2dxds&\leq \int_1^t\int_{\mathbb R^3}|\nabla^2 d|^4+|u|^4dxds\\&\leq \bar E
+\int_1^t(\int_{\mathbb R^3}|u|^2dx)^\frac{1}{2}(\int_{\mathbb
R^3}|\nabla u|^2dx)^\frac{3}{2}ds\\&\leq M(\bar
E+C_0^\frac{3}{2}\bar A^\frac{1}{2}).\end {array}$$ Taking the above
results into (3.11), then (3.7) follows. Thus the proof of lemma is completed.\\

Next we derive preliminary bounds for $\dot u$ and $\nabla d_t$ in
$L^\infty([0,T];L^2(\mathbb R^3))$.\\
\\{\bf Lemma 3.4.} Assume that the hypotheses and notations of
Proposition 3.1 are in force. Then for $0<t\leq 1\wedge T$,
$$\begin {array}{rl}&\sup\limits_{0< s\leq t}\sigma^5\int_{\mathbb R^3} (|\dot u|^2+|\nabla d_t|^2)dx+
\int_0^t\int_{\mathbb R^3}\sigma^5(|\nabla\dot u|^2+|\nabla^2
d_t|^2)dxds\\\leq &M[C_0+E+C_0^{\frac{q-4}{q-2}}B_q^{\frac{2}{q-2}}+
C_0^{\frac{q-6}{q-2}}B_q^{\frac{4}{q-2}}+C_0^{\frac{q-4}{q-2}}B_q^{\frac{2}{q-2}}(E+C_0^{\frac{q-4}{q-2}}B_q^{\frac{2}{q-2}})\\&+
(C_0^{\frac{q-6}{q-2}}B_q^{\frac{4}{q-2}})^\frac{1}{3}A],\end
{array}\eqno(3.13)$$ and if $T>1$ and $1\leq t\leq T$, then
$$\begin {array}{rl}&\sup\limits_{1\leq s\leq t}\int_{\mathbb R^3} (|\dot u|^2+|\nabla d_t|^2)dx+
\int_0^t\int_{\mathbb R^3}(|\nabla\dot u|^2+|\nabla^2
d_t|^2)dxds\\\leq &M\{C_0+C_0\bar A\bar
E+C_0^\frac{2(q-3)}{3(q-2)}\bar B_q^\frac{2}{3(q-2)}\bar A+\bar
E\}+A(1).\end {array}\eqno(3.14)$$
\\{\bf Proof.} By the definition of material derivative, we can
write (1.2) as follows,
$$\rho\dot u+\nabla(P(\rho))=\mu\triangle u+\lambda\nabla div u-\nabla d\cdot\triangle d.\eqno(3.15)$$
Differentiation (3.15) with respect to $t$ and using (1.1), we have
$$\begin {array}{rl}&\rho\dot u_t+\rho u\cdot\nabla\dot u+\nabla(P(\rho)_t)+(\nabla d\cdot\triangle d)_t\\
=&\mu\triangle \dot u+\lambda\nabla div \dot u-[\mu\triangle
(u\cdot\nabla u)+\lambda\nabla div(u\cdot\nabla u)]\\&
+div[(\mu\triangle u+\lambda\nabla div u)\otimes u-\nabla
P(\rho)\otimes u-(\nabla d\cdot\triangle d)\otimes u].\end
{array}\eqno(3.16)$$ Multiplying (3.16) by $\sigma^5\dot u$, and
integrating over $\mathbb R^3\times [0,t]$, we obtain that for
$0\leq t\leq 1\wedge T$,
$$\begin {array}{rl}&\sup\limits_{0<s\leq t}\sigma^5\int_{\mathbb R^3}|\dot u|^2dx+\int_0^t\int_{\mathbb R^3}\sigma^5|\nabla\dot u|^2dxds
\\ \leq &M[C_0+E+C_0^{\frac{q-4}{q-2}}B_q^{\frac{2}{q-2}}+C_0^{\frac{q-6}{q-2}}B_q^{\frac{4}{q-2}}
\\&+\int_0^t\int_{\mathbb R^3}\sigma^5|\nabla d|^2(|\nabla
d_t|^2+|u|^2|\nabla^2d|^2)dxds].\end {array}\eqno(3.17)$$ Next we
differentiate (3.9) with respect to $t$, multiply by $\sigma^5\nabla
d_t$ and integrate over $\mathbb R^3\times [0,t]$ to obtain
$$\begin {array}{rl}&\frac{1}{2}\sigma^5\int_{\mathbb R^3}|\nabla d_t|^2dx+\int_0^t\int_{\mathbb R^3}\sigma^5|\nabla^2d_t|^2dxds\\
=&\frac{5}{2}\int_0^t\int_{\mathbb R^3}\sigma^4\sigma'|\nabla
d_t|^2dxds+\int_0^t\int_{\mathbb R^3}\sigma^5\nabla(|\nabla
d|^2d-u\cdot\nabla d)_t\nabla d_tdxds.\end {array}$$ Adding this to
(3.17), integrating by parts, using Cauchy's inequality, we then
have
$$\begin {array}{rl}&\sup\limits_{0<s\leq t}\sigma^5\int_{\mathbb R^3}|\dot u|^2+|\nabla d|^2dx
+\int_0^t\int_{\mathbb R^3}\sigma^5(|\nabla\dot
u|^2+|\nabla^2d|^2)dxds
\\ \leq &M[C_0+E+C_0^{\frac{q-4}{q-2}}B_q^{\frac{2}{q-2}}+C_0^{\frac{q-6}{q-2}}B_q^{\frac{4}{q-2}}
\\&+\int_0^t\int_{\mathbb R^3}\sigma^5|\nabla d|^2|u|^2(|\nabla u|^2+|\nabla^2d|^2)+\sigma^5|\nabla d|^4|d_t|^2 dxds
\\&+\int_0^t\int_{\mathbb R^3}\sigma^5(|\nabla d|^2|\dot u|^2+|\nabla d|^2|\nabla d_t|^2+|\nabla d_t|^2|u|^2)dxds].\end {array}\eqno(3.18)$$
By (2.2), the terms on right side can be bounded by
$$\begin {array} {rl}&\int_0^t\int_{\mathbb R^3}\sigma^5|\nabla d|^2|u|^2|\nabla
u|^2dxds\\ \leq &\int_0^t\int_{\mathbb R^3}\sigma^5(|\nabla
d|^8+|u|^8+|\nabla u|^4)dxds\\ \leq& E+\sup\limits_{0\leq s\leq
t}\|(u,\nabla d)(s)\|^4_{L^4(\mathbb
R^3)}\int^t_0\sigma^5\|(u,\nabla d)(s)\|^4_{L^\infty(\mathbb R^3)}ds\\
\leq
&E+M[C_0^{\frac{q-4}{q-2}}B_q^{\frac{2}{q-2}}][E+C_0^{\frac{q-4}{q-2}}B_q^{\frac{2}{q-2}}].
\end {array}$$
By (1.1c) and (3.12), we have
$$\begin {array} {rl}&\int_0^t\int_{\mathbb R^3}\sigma^5|\nabla d|^4|d_t|^2dxds\\ \leq &M\int_0^t\int_{\mathbb R^3}\sigma^5(|\nabla
d|^4|\nabla^2d|^2+|\nabla d|^4|u|^2|\nabla d|^2+|\nabla d|^8)dxds\\
\leq&M\int_0^t\int_{\mathbb R^3}\sigma^5(|\nabla
d|^4|\nabla^2d|^2+|\nabla^2 d|^2|u|^2|\nabla d|^2)dxds\\ \leq
&ME+M[C_0^{\frac{q-4}{q-2}}B_q^{\frac{2}{q-2}}][E+C_0^{\frac{q-4}{q-2}}B_q^{\frac{2}{q-2}}].
\end {array}$$
The last term on the right side in (3.18) can be bounded by
$$\begin {array} {rl}&\int_0^t\int_{\mathbb R^3}\sigma^5|\nabla d|^2|\dot u|^2dxds\\ \leq &
(\int_0^t\int_{\mathbb R^3}|\nabla
d|^6dxds)^\frac{1}{3}(\int_0^t\int_{\mathbb
R^3}\sigma^\frac{15}{2}|\dot u|^3dxds)^\frac{2}{3}\\ \leq &
(C_0^{\frac{q-6}{q-2}}B_q^{\frac{4}{q-2}})^\frac{1}{3}
(\int_0^t\sigma^\frac{15}{2}\|\dot u\|_{L^2(\mathbb
R^3)}^\frac{3}{2}\|\nabla\dot u\|_{L^2(\mathbb
R^3)}^\frac{3}{2}dxds)^\frac{2}{3}\\ \leq &
(C_0^{\frac{q-6}{q-2}}B_q^{\frac{4}{q-2}})^\frac{1}{3}
(\int_0^t\sigma^{15}\|\dot u\|_{L^2(\mathbb R^3)}^6dxds)^\frac{1}{6}
(\int_0^t\sigma^5\|\nabla\dot u\|_{L^2(\mathbb
R^3)}^2dxds)^\frac{1}{2}\\ \leq &
(C_0^{\frac{q-6}{q-2}}B_q^{\frac{4}{q-2}})^\frac{1}{3}A.
\end {array}$$
The other integrals on the right side of (3.18) are bounded in a
similar way, and (3.13) follows.

For $1\leq t\leq T$, as in (3.18), we have
$$\begin {array}{rl}&\sup\limits_{1\leq s\leq t}\int_{\mathbb R^3}|\dot u|^2+|\nabla d|^2dx
+\int_1^t\int_{\mathbb R^3}(|\nabla\dot u|^2+|\nabla^2d|^2)dxds\\
\leq &M[\int_1^t\int_{\mathbb R^3}|\nabla d|^2|u|^2(|\nabla
u|^2+|\nabla^2d|^2)+|\nabla d|^4|d_t|^2 dxds
\\&+\int_1^t\int_{\mathbb R^3}(|\nabla d|^2|\dot u|^2+|\nabla d|^2|\nabla d_t|^2+|\nabla d_t|^2|u|^2)dxds
\\ &+C_0+\bar E]+A(1).\end {array}\eqno(3.19)$$
The terms on the right can be bounded by
$$\begin {array} {rl}&\int_1^t\int_{\mathbb R^3}|\nabla d|^2|u|^2|\nabla
u|^2dxds\\ \leq &\int_1^t\int_{\mathbb R^3}(|\nabla
d|^8+|u|^8+|\nabla u|^4)dxds\\ \leq& \bar E+\sup\limits_{1\leq s\leq
t}\|(u,\nabla d)(s)\|^\frac{10}{3}_{L^\frac{10}{3}(\mathbb
R^3)}\int^t_1\|(u,\nabla d)(s)\|^\frac{14}{3}_{L^\infty(\mathbb R^3)}ds\\
\leq & \bar E+\sup\limits_{1\leq s\leq t}[\|(u,\nabla
d)(s)\|^\frac{4}{3}_{L^2(\mathbb R^3)}\|\nabla (u,\nabla
d)(s)\|^2_{L^2(\mathbb R^3)}]\\&\times\int^t_1\|(u,\nabla
d)(s)\|^\frac{2}{3}_{L^2 (\mathbb R^3)}\|\nabla(u,\nabla
d)(s)\|^4_{L^4 (\mathbb R^3)}ds\\ \leq & \bar E+MC_0\bar A\bar E.
\end {array}$$
By (1.1c) and (3.12), we have
$$\begin {array} {rl}&\int_1^t\int_{\mathbb R^3}|\nabla d|^4|d_t|^2dxds\\ \leq &M\int_1^t\int_{\mathbb R^3}(|\nabla
d|^4|\nabla^2d|^2+|\nabla d|^4|u|^2|\nabla d|^2+|\nabla d|^8)dxds\\
\leq&M\int_1^t\int_{\mathbb R^3}(|\nabla d|^4|\nabla^2d|^2+|\nabla^2
d|^2|u|^2|\nabla d|^2)dxds\\ \leq &M\bar E+MC_0\bar A\bar E.
\end {array}$$
The last term on the right side in (3.19) can be bounded by
$$\begin {array} {rl}&\int_1^t\int_{\mathbb R^3}|\nabla d|^2|\dot u|^2dxds\\ \leq &
\int_1^t(\int_{\mathbb R^3}|\nabla
d|^3dx)^\frac{2}{3}(\int_{\mathbb R^3}|\dot u|^6dx)^\frac{1}{3}ds\\
\leq &\sup\limits_{1\leq s\leq t} \int_{\mathbb R^3}|\nabla
d|^3dx)^\frac{2}{3}\int_1^t\int_{\mathbb R^3}|\nabla \dot u|^2dxds
\\ \leq &C_0^\frac{2(q-3)}{3(q-2)}\bar B_q^\frac{2}{3(q-2)}\bar A.
\end {array}$$
The other integrals on the right side of (3.19) are bounded in a
similar way, and (3.14) follows. Thus the proof of lemma is
completed.\\

Next we derive a number of auxiliary estimates needed to close the
bounds in the previous two lemmas. We begin with a bound for the
vorticity $\omega$.\\
\\{\bf Lemma 3.5.} Assume that the hypotheses and notations of
Proposition 3.1 are in force. Then for $0<t\leq 1\wedge T$,$$\begin
{array}{rl}&\sup\limits_{0<s\leq t}\int_{\mathbb
R^3}\sigma^5(|\nabla F|^2+|\nabla\omega|^2)dx+\int_0^t\int_{\mathbb
R^3}\sigma(|\nabla F|^2+|\nabla\omega|^2)dxds\\ \leq
&M(D+E^\frac{2}{3}C_0^\frac{q-6}{3q-6}B_q^\frac{4}{3q-4}+\sup\limits_{0<s\leq
t}\int_{\mathbb R^3}\sigma^5|\dot u|^2dx+\int_0^t\int_{\mathbb
R^3}\sigma|\dot u|^2dxds),
\end {array}\eqno (3.20)$$
and if $T>1$ and $1\leq t\leq T$, then
$$\begin
{array}{rl}&\sup\limits_{1\leq s\leq t}\int_{\mathbb R^3}(|\nabla
F|^2+|\nabla\omega|^2)dx+\int_1^t\int_{\mathbb R^3}(|\nabla
F|^2+|\nabla\omega|^2)dxds\\ \leq &M(\bar D+\bar
E+\sup\limits_{1\leq s\leq t}\int_{\mathbb R^3}|\dot
u|^2dx+\int_1^t\int_{\mathbb R^3}|\dot u|^2dxds).
\end {array}\eqno (3.21)$$
\\
{\bf Proof.} By (2.10) and the definition of $D,E$, we can easily
get (3.20) and (3.21). The proof of lemma is completed.\\

Next we derive an estimate for the functional $B_q$.\\
\\{\bf Lemma 3.6.} Assume that the hypotheses and notations of
Proposition 3.1 are in force. Then for any $0<t\leq T,$
$$B_q\leq M(C_0^\frac{p-q}{p-2}N^\frac{q-2}{p-2}+C_0^\frac{q+3}{3(q-2)}B_q^\frac{3q-11}{3(q-2)}+C_0^\frac{q-3}{3(q-2)}B_q^\frac{3q-5}{3(q-2)}).\eqno (3.22)$$
\\{\bf Proof.} We multiply (1.1b) by $|u|^{q-2}u$ and integrate over
$\mathbb R^3\times (0,t)$ to obtain that
$$\begin {array}{rl}&q^{-1}\int_{\mathbb R^3}\rho |u|^qdx|_0^t+\int_0^t\int_{\mathbb R^3}\mu |u|^{q-2}|\nabla
u|^2dxds
\\&+\int_0^t\int_{\mathbb R^3}[\frac{1}{4}\mu(q-2) |u|^{q-4}|\nabla
(|u|^2)|^2+\lambda |u|^{q-2}(div u)^2]dxds\\=& \int_0^t\int_{\mathbb
R^3}[(P(\rho)-\tilde P(\rho))div(|u|^{q-2}u)-|u|^{q-2}u\nabla
d\cdot\triangle d]dxds\\&-\int_0^t\int_{\mathbb
R^3}\frac{1}{2}\lambda(q-2) |u|^{q-4}(div u)u\cdot\nabla(|u|^2)dxds.
\end {array}\eqno (3.23)$$
For any $\eta>0$,
$$\begin {array}{rl}&|-\int_0^t\int_{\mathbb
R^3}\frac{1}{2}\lambda(q-2) |u|^{q-4}(div u)u\cdot\nabla(|u|^2)dxds|
\\ \leq&\frac{1}{2}\lambda(q-2)\int_0^t\int_{\mathbb
R^3} |u|^\frac{q-2}{2}|(div u)||u|^\frac{q-4}{2}|\nabla(|u|^2)|dxds|
\\ \leq&\frac{1}{4}\lambda(q-2)[\eta\int_0^t\int_{\mathbb
R^3} |u|^{q-2}|(div u)|^2dxds+\eta^{-1}\int_0^t\int_{\mathbb
R^3}|u|^{q-4}|\nabla(|u|^2)|^2dxds],
\end {array}$$
so if we choose
$$\frac{1}{4}\lambda(q-2)\eta=\beta\mu+\lambda$$
for a positive $\beta$ to be determined, then the term in question
will be bounded by
$$\begin {array}{rl}&3\beta\mu\int_0^t\int_{\mathbb
R^3} |u|^{q-2}|\nabla u|^2dxds+\lambda\int_0^t\int_{\mathbb R^3}
|u|^{q-2}|div
u|^2dxds\\&+\frac{[\frac{1}{4}\lambda(q-1)]^2}{\beta\mu+\lambda}\int_0^t\int_{\mathbb
R^3}|u|^{q-4}|\nabla(|u|^2)|^2dxds.
\end {array}$$
Substituting this into (3.23), we then get
$$\begin {array}{rl}&q^{-1}\int_{\mathbb R^3}\rho |u|^qdx|_0^t+\mu(1-3\beta)\int_0^t\int_{\mathbb R^3} |u|^{q-2}|\nabla
u|^2dxds
\\&+[\frac{1}{4}\mu(q-2)+\frac{[\frac{1}{4}\lambda(q-1)]^2}{\beta\mu+\lambda}]\int_0^t\int_{\mathbb R^3} |u|^{q-4}|\nabla
(|u|^2)|^2dxds\\ \leq& |\int_0^t\int_{\mathbb R^3}(P(\rho)-\tilde
P(\rho))div(|u|^{q-2}u)dxds|\\&+|\int_0^t\int_{\mathbb
R^3}|u|^{q-2}u\nabla d\cdot\triangle ddxds|.
\end {array}$$
Recall that $q\in[6,p)$, thus (1.9) holds with $p$ replaced by $q$,
and this is the condition that brackets on the left here is positive
when $\beta=\frac{1}{3}$. It follows this term is positive for some
$\beta\in (0,\frac{1}{3})$, which we now fix. It then follows that
$$\begin {array}{rl}&q^{-1}\int_{\mathbb R^3}\rho |u|^qdx|_0^t+\int_0^t\int_{\mathbb R^3} |u|^{q-2}|\nabla
u|^2dxds +\int_0^t\int_{\mathbb R^3} |u|^{q-4}|\nabla
(|u|^2)|^2dxds\\ \leq&M[ |\int_0^t\int_{\mathbb R^3}(P(\rho)-\tilde
P(\rho))div(|u|^{q-2}u)dxds|+|\int_0^t\int_{\mathbb
R^3}|u|^{q-2}u\nabla d\cdot\triangle ddxds|].
\end {array}\eqno(3.24)$$
We multiply (3.9) by $|\nabla d|^{q-2}\nabla d$ and integrate over
$\mathbb R^3\times(0,t)$ to obtain that
$$\begin {array}{rl}&q^{-1}\int_{\mathbb R^3} |\nabla d|^q dx|_0^t+\int_0^t\int_{\mathbb R^3} |\nabla
d|^{q-2}|\nabla^2 d|^2dxds \\&+\int_0^t\int_{\mathbb R^3}(q-2)
|\nabla d|^{q-4}|\nabla (|\nabla d|^2)|^2dxds\\ =&
\int_0^t\int_{\mathbb R^3}|\nabla d|^{q-2}\nabla d\nabla (|\nabla
d|^2d-u\cdot\nabla d)dxds.
\end {array}\eqno(3.25)$$
Adding (3.25) to (3.24) and applying the Cauchy's inequality in an
elementary way we then obtain
$$\begin {array}{rl}&\int_{\mathbb R^3} |u|^q+|\nabla d|^qdx+\int_0^t\int_{\mathbb R^3} |u|^{q-2}|\nabla
u|^2dxds +\int_0^t\int_{\mathbb R^3} |u|^{q-4}|\nabla
(|u|^2)|^2dxds\\&+\int_0^t\int_{\mathbb R^3} |\nabla
d|^{q-2}|\nabla^2 d|^2dxds +\int_0^t\int_{\mathbb R^3} |\nabla
d|^{q-4}|\nabla (|\nabla d|^2)|^2dxds\\ \leq&M[\int_{\mathbb R^3}
|u_0|^q+|\nabla d_0|^qdx+ |\int_0^t\int_{\mathbb R^3}(P(\rho)-\tilde
P(\rho))div(|u|^{q-2}u)dxds|\\&+|\int_0^t\int_{\mathbb
R^3}|u|^{q-2}u\nabla d\cdot\triangle ddxds|+|\int_0^t\int_{\mathbb
R^3}|\nabla d|^{q-2}\nabla d\nabla (|\nabla d|^2d-u\cdot\nabla
d)dxds|]\\=&\sum^4_{i=1}I_i.
\end {array}\eqno(3.26)$$
Since $q\in [6,\min\{p,12\})$, then by H\"older's inequality and
Sobolev's inequality, we have
$$I_1\leq (\int_{\mathbb R^3}
|u_0|^2+|\nabla d_0|^2dx)^\frac{p-q}{p-2}(\int_{\mathbb R^3}
|u_0|^p+|\nabla d_0|^pdx)^\frac{q-2}{p-2}\leq
MC_0^\frac{p-q}{p-2}N^\frac{q-2}{p-2},$$
and
$$\begin {array}{rl}I_2\leq &[\int_0^t\int_{\mathbb R^3} |u|^{2q-4}dxds]^\frac{1}{2}[\int_0^t\int_{\mathbb
R^3} |\nabla u|^2dxds]^\frac{1}{2}\\ \leq
&C_0^\frac{1}{2}[\int_0^t(\int_{\mathbb R^3}
|u|^{3q}dx)^\frac{1}{3}(\int_{\mathbb R^3}
|u|^{\frac{3}{2}(q-4)}dx)^\frac{2}{3}ds]^\frac{1}{2}\\ \leq&
C_0^\frac{1}{2}[\int_0^t(\int_{\mathbb R^3} |u|^{q-2}|\nabla
u|^2dx)(\int_{\mathbb R^3}
|u|^{\frac{3}{2}(q-4)}dx)^\frac{2}{3}ds]^\frac{1}{2}\\ \leq&
C_0^\frac{1}{2}B_q^\frac{1}{2}\sup\limits_{0<s\leq t}[\int_{\mathbb
R^3} |u(s)|^2dx]^\frac{12-q}{6(q-2)}\int_{\mathbb R^3}
|u(s)|^qdx]^\frac{3q-16}{6(q-2)}\\
\leq&MC_0^\frac{q+3}{3(q-2)}B_q^\frac{3q-11}{3(q-2)}.
\end
{array}$$
Using the fact
$$\nabla d\cdot\triangle d=div(\nabla d\otimes\nabla d)-\frac{1}{2}\nabla d|\nabla d|^2$$
and integrating by parts, we have
$$\begin {array}{rl}I_3\leq & [\int_0^t\int_{\mathbb R^3} (|u|^{q-2}|\nabla
u|^2+|u|^{q-4}|\nabla(|u|^2)|^2)dxds]^\frac{1}{2}[\int_0^t\int_{\mathbb
R^3} |\nabla d|^4|u|^{q-2}dxds]^\frac{1}{2}\\ \leq
&MB_q^\frac{1}{2}[(\int_0^t\int_{\mathbb R^3}|\nabla d|^{q+2}dxds
)^\frac{1}{2}+(\int_0^t\int_{\mathbb R^3}| u|^{q+2}dxds
)^\frac{1}{2}]\\ \leq &MB_q^\frac{1}{2}[(\int_0^t(\int_{\mathbb
R^3}|\nabla d|^{3q}dx)^\frac{1}{3}(\int_{\mathbb R^3}|\nabla
d|^{3}dx)^\frac{2}{3}ds )^\frac{1}{2}\\&+(\int_0^t(\int_{\mathbb
R^3}|u|^{3q}dx)^\frac{1}{3}(\int_{\mathbb
R^3}|u|^{3}dx)^\frac{2}{3}ds )^\frac{1}{2}]\\ \leq
&MC_0^\frac{q-3}{3(q-2)}B_q^\frac{3q-5}{3(q-2)}.
\end {array}$$
Similarly, we have
$$\begin {array}{rl}I_4\leq & M[\int_0^t\int_{\mathbb R^3} |\nabla^2 d||\nabla d|^{q-2}(|\nabla
d|^2+|u||\nabla d|^2)dxds]\\ \leq &M[\int_0^t\int_{\mathbb R^3}
(|\nabla d|^{q-2}|\nabla^2
d|^2dxds]^\frac{1}{2}[\int_0^t\int_{\mathbb R^3} |\nabla
d|^{q+2}+|\nabla d|^{q}|u|^2dxds]^\frac{1}{2} \\ \leq
&MC_0^\frac{q-3}{3(q-2)}B_q^\frac{3q-5}{3(q-2)}.
\end {array}$$
Substituting these results into (3.26) gives (3.22). Thus the proof of lemma is completed.\\

Next we derive a bound for the functional $D$ and $\bar D$.\\
\\{\bf Lemma 3.7.} Assume that the hypotheses and notations of
Proposition 3.1 are in force. Then for $0<t\leq 1\wedge T$,
$$D\leq M[C_0^\frac{q-4}{2q-4}B_q^\frac{1}{q-2}A+A^2+A^5],\eqno(3.27)$$
and if $T>1$ and $1\leq t\leq T$, then
$$\bar D\leq M[C_0^\frac{1}{4}\bar A^\frac{7}{4}+\bar A^2+\bar
A^5].\eqno(3.28)$$\\
{\bf Proof.} We give the proof of (3.27), that of (3.28) being
similar. By Lemma 2.1, we have
$$\begin {array}{rl}&\int_{\mathbb R^3}\sigma^5|\nabla^2 d|^2|\nabla d|^2dx\\ \leq&\sigma^4
\|\nabla
d(\cdot,t)\|^2_{\infty}[\sigma\|\nabla^2d(\cdot,t)\|^2_{L^2(\mathbb
R^3)}]\\ \leq& MA[\sigma^4\|\nabla d(\cdot,t)\|^2_{L^4(\mathbb
R^3)}+\sigma^4\|\nabla^2d(\cdot,t)\|^2_{L^4(\mathbb
R^3)}]\\
\leq
&MA[C_0^\frac{q-4}{2q-4}B_q^\frac{1}{q-2}+\|\sigma^\frac{1}{2}\nabla^2d(\cdot,t)\|^\frac{1}{2}_{L^2(\mathbb
R^3)}\|\sigma^\frac{5}{2}\nabla^3d(\cdot,t)\|^\frac{3}{2}_{L^2(\mathbb
R^3)}]\\
\leq
&MA[C_0^\frac{q-4}{2q-4}B_q^\frac{1}{q-2}+A^\frac{1}{4}\|\sigma^\frac{5}{2}\nabla^3d(\cdot,t)\|^\frac{3}{2}_{L^2(\mathbb
R^3)}].
\end {array}$$
From (3.9) and (3.12), we have
$$\begin {array}{rl}\|\sigma^\frac{5}{2}\nabla^3d(\cdot,t)\|^2_{L^2(\mathbb R^3)}&\leq M\int_{\mathbb R^3}\sigma^5(|\nabla d_t|^2+|u|^2|\nabla^2d|^2+
|\nabla u|^2|\nabla d|^2+|\nabla^2d|^2|\nabla d|^2)dx\\&\leq
M(A+D),\end {array}$$ so that
$$\int_{\mathbb R^3}\sigma^5|\nabla^2 d|^2|\nabla d|^2dx\leq MA[C_0^\frac{q-4}{2q-4}B_q^\frac{1}{q-2}+A^\frac{1}{4}(A+D)^\frac{3}{4}].$$
The other terms included in $D$ are estimated in exactly the same
way, and (3.27) follows. The proof of lemma is completed.\\

The following lemma contains the required bound for the pressure
term in (2.11), which has been proved in Hoff [8, Lemma 3.3].\\
\\{\bf Lemma 3.8.} Assume that the hypotheses and notations of
Proposition 3.1 are in force. Then it holds
$$\int_0^t\int_{\mathbb R^3} \sigma^5|\rho-\tilde\rho|^4dxds\leq M[C_0+\int_0^t\int_{\mathbb R^3} \sigma^5|F|^4dxds.\eqno (3.29)$$

We can now obtain the required estimates for the functional $E$ and
$\bar E$.\\
\\{\bf Lemma 3.9.} Assume that the hypotheses and notations of
Proposition 3.1 are in force. Then there are polynomials $\varphi_1$
and $\varphi_2$ whose degrees and coefficients depend on the same
$M$ quantities as $M$ in the statement of Proposition 3.1 such that:
for $0<t\leq 1\wedge T$ $$E\leq
M[\varphi_1(C_0)+\varphi_2(A+B_q)],\eqno(3.30)$$ and if $T>1$ and
$1\leq t\leq T$, then
$$\bar E\leq M[\varphi_1(C_0+A(1)+B_q(1))+\varphi_2(\bar A+B_q)].\eqno(3.31)$$
The polynomial $\varphi_1$ contains no constant term and the
monomials
in $\varphi_2$ all have degrees strictly greater than 1.\\
\\{\bf Proof.} Since the term $|\sum\limits_{1\leq k_i,j_m\leq 3}\int_0^t\int_{\mathbb R^3}\sigma u^{j_1}_{x_{k_1}}u^{j_2}_{x_{k_2}}u^{j_3}_{x_{k_3}}dxds|$
has been bounded exactly in Hoff [8]. So here we just bound the
other terms for simplicity.

First for $0< t\leq 1\wedge T$, from (2.1), (3.9) and (3.12) we have
$$\begin {array} {rl}&
\int_0^t\int_{\mathbb R^3}\sigma^\frac{3}{2}|\nabla^2d|^3dxds\\ \leq
&M (\int_0^t\int_{\mathbb
R^3}\sigma|\nabla^2d|dxds)^\frac{3}{4}(\int_0^t\int_{\mathbb
R^3}\sigma|\nabla^3d|^2dxds)^\frac{3}{4}\\ \leq &MA^\frac{3}{4}
(\int_0^t\int_{\mathbb R^3}\sigma(|\nabla d_t|^2+|\nabla
d|^2|\nabla^2d|^2+|u|^2|\nabla^2d|^2+|\nabla d|^2|\nabla u|^2)dxds)^\frac{3}{4}\\
\leq &MA^\frac{3}{4}(A+ \int_0^t(\int_{\mathbb R^3}|\nabla
d|^6+|u|^6+\sigma^\frac{3}{2}(|\nabla^2d|^3+|\nabla
u|^3)dxds)^\frac{3}{4}\\
\leq&MA(A+C_0^\frac{q-6}{q-2}B_q^\frac{4}{q-2}+E)^\frac{3}{4},
\end
{array}$$ and
$$\begin {array} {rl}&
\int_0^t\int_{\mathbb R^3}\sigma^5|\nabla^2d|^4dxds\\ \leq &M
(\int_0^t\sigma^5(\int_{\mathbb
R^3}|\nabla^2d|dx)^\frac{1}{2}(\int_{\mathbb
R^3}|\nabla^3d|^2dx)^\frac{3}{2}ds\\ \leq &M
(\int_0^t(\sigma\int_{\mathbb
R^3}|\nabla^2d|dx)^\frac{1}{2}(\sigma\int_{\mathbb
R^3}|\nabla^3d|^2dx)^\frac{3}{4}(\sigma^5\int_{\mathbb
R^3}|\nabla^3d|^2dx)^\frac{3}{4}ds
\\
\leq&MA(A+C_0^\frac{q-6}{q-2}B_q^\frac{4}{q-2}+E)^\frac{3}{4}(A+D)^\frac{3}{4}.
\end
{array}$$ From Lemma 2.4, Lemma 3.2, Lemma 3.5 and the definition of
$F,\omega$, we have
$$\begin {array}{rl}&\int_0^t\int_{\mathbb R^3} \sigma^5|\nabla u|^4dxds\\ \leq &M[\int_0^t\int_{\mathbb R^3} \sigma^5(|\rho-\tilde\rho|^4+|F|^4+|\omega|^4)dxds]
\\\leq & M[C_0+(\sup\limits_{0\leq s\leq t}\int_{\mathbb R^3} \sigma(|F|^2+|\omega|^2)dx
\int_{\mathbb R^3} \sigma^5(|\nabla
F|^2+|\nabla\omega|^2)dx)^\frac{1}{2}\\&\times(\int_0^t\int_{\mathbb
R^3} \sigma(|\nabla F|^2+|\nabla\omega|^2)dx)\\ \leq &M
[C_0+(C_0+A)^\frac{1}{2}(A+D)^\frac{1}{2}(A+E^\frac{2}{3}C_0^\frac{q-6}{3q-6}B_q^\frac{4}{3q-4})]
,\end {array}\eqno $$
$$\begin {array}{rl}&\int_0^t\int_{\mathbb R^3} \sigma^\frac{3}{2}|\nabla u|^3dxds\\ \leq &M[\int_0^t\int_{\mathbb R^3} \sigma^\frac{3}{2}(|\rho-\tilde\rho|^3+|F|^3+|\omega|^3)dxds]
\\\leq & M[C_0+\int_0^t\sigma^{\frac{3}{2}}(\int_{\mathbb
R^3}(|F|^2+|\omega|^2)dx)^\frac{3}{4} (\int_{\mathbb R^3} (|\nabla
F|^2+|\nabla\omega|^2)dx)^\frac{3}{4}ds\\ \leq &M
[C_0+(C_0+A)^\frac{3}{4}(A+E^\frac{2}{3}C_0^\frac{q-6}{3q-6}B_q^\frac{4}{3q-4})^\frac{3}{4}]
.\end {array}$$ Thus combining the above results and Lemma (3.7), we
yields (3.30).

Now for $1\leq t\leq T, $ if we take $q=4$ in (3.25) and integrate
by parts to obtain that
$$\begin {array}{rl}&\int_{\mathbb R^3} |\nabla d|^4 dx+\int_1^t\int_{\mathbb R^3} |\nabla
d|^{2}|\nabla^2 d|^2dxds +\int_1^t\int_{\mathbb R^3} |\nabla
(|\nabla d|^2)|^2dxds\\ \leq & M[\int_{\mathbb R^3} |\nabla
d(\cdot,1)|^4 dx+|\int_1^t\int_{\mathbb R^3}\triangle d|\nabla
d|^{4}ddxds|+\int_1^t\int_{\mathbb R^3}|\nabla d|^{3}|\nabla^2
d||u|dxds]\\ \leq&
M[C_0^\frac{q-4}{q-2}B_q(1)^\frac{2}{q-2}+\int_1^t\int_{\mathbb
R^3}|\triangle d+|\nabla d|^2d|^2|\nabla
d|^{2}dxds\\&+(\int_0^t\int_{\mathbb R^3}|\nabla
d|^\frac{18}{5}|\nabla^2
d|^\frac{6}{5}dxds)^\frac{5}{6}(\int_0^t\int_{\mathbb
R^3}|u|^6dxds)^\frac{1}{6}.\\ \leq&
M[C_0^\frac{q-4}{q-2}B_q(1)^\frac{2}{q-2}+\sup\limits_{1\leq s\leq
t}\|\nabla d\|_{L^\infty(\mathbb R^3)}^{2}\int_1^t\int_{\mathbb
R^3}|\triangle d+|\nabla d|^2d|^2dxds\\&+(\int_0^t\int_{\mathbb
R^3}|\nabla d|^2|\nabla^2
d|^2dxds)^\frac{5}{6}(\int_0^t\int_{\mathbb R^3}|\nabla
u|^2dxds)^\frac{1}{2}\\ \leq
&M[C_0^\frac{q-4}{q-2}B_q(1)^\frac{2}{q-2}+C_0(C_0^\frac{q-4}{2q-4}B_q^\frac{1}{q-2}+\bar
A^\frac{1}{4}(\bar A+\bar D)^\frac{3}{4})+C_0^\frac{1}{2}\bar
E^\frac{5}{6}].
\end {array}$$
Multiplying (3.9) by $\nabla\triangle d$ and integrating over
$\mathbb R^3$, we have
$$\begin {array} {rl}
\int_{\mathbb R^3}|\nabla^3d|^2dx\leq&\int_{\mathbb R^3}|\nabla
d_t||\nabla^3d|dx+2\int_{\mathbb R^3}|\nabla
d||\nabla^2d||\nabla^3d|dx\\&+\int_{\mathbb R^3}\nabla u\cdot\nabla
d\nabla\triangle d+u\cdot\nabla \nabla d\nabla\triangle ddx\\
=&\int_{\mathbb R^3}|\nabla d_t||\nabla^3d|dx+2\int_{\mathbb
R^3}|\nabla d||\nabla^2d||\nabla^3d|dx\\&+\int_{\mathbb R^3}\nabla
u\cdot\nabla d\nabla\triangle d-\nabla u\cdot\nabla \nabla
d\triangle d+\frac{1}{2}(divu)|\triangle d|^2dx.
\end
{array}$$ By Cauchy's inequality, we have
$$\begin {array} {rl}
\int_{\mathbb R^3}|\nabla^3d|^2dx\leq&\int_{\mathbb R^3}|\nabla
d_t|^2+|\nabla d|^2|\nabla^2d|^2+|\nabla d|^2|\nabla u|^2+|\nabla
u||\nabla^2 d|^2dx.
\end
{array}$$ Thus we have
$$\begin {array} {rl}&
\int_1^t\int_{\mathbb R^3}|\nabla^2d|^4dxds\\ \leq &M
(\int_1^t(\int_{\mathbb R^3}|\nabla^2d|dx)^\frac{1}{2}(\int_{\mathbb
R^3}|\nabla^3d|^2dx)^\frac{3}{2}ds\\ \leq &M\bar
A\int_1^t(\int_{\mathbb R^3}|\nabla d_t|^2+|\nabla
d|^2|\nabla^2d|^2+|\nabla d|^2|\nabla u|^2+|\nabla u||\nabla^2
d|^2dx)\\&\times (\int_{\mathbb R^3}|\nabla d_t|^2+|\nabla
d|^2|\nabla^2d|^2+|\nabla d|^2|\nabla u|^2+| u|^2|\nabla^2
d|^2dx)^\frac{1}{2}ds
\\
\leq&M\bar A(\bar A+\bar D)^\frac{1}{2}(\int_1^t(\int_{\mathbb
R^3}|\nabla d_t|^2+|\nabla
d|^2|\nabla^2d|^2dxds\\&+\int_1^t(\int_{\mathbb R^3}|\nabla
d|^2|\nabla u|^2+|\nabla u||\nabla^2 d|^2dxds)\\ \leq & M\bar A(\bar
A+\bar
D)^\frac{1}{2}(A+C_0^\frac{q-4}{q-2}B_q(1)^\frac{2}{q-2}\\&+C_0(C_0^\frac{q-4}{2q-4}B_q^\frac{1}{q-2}+\bar
A^\frac{1}{4}(\bar A+\bar D)^\frac{3}{4})+C_0^\frac{1}{2}\bar
E^\frac{5}{6}+C_0^\frac{1}{2}\bar E^\frac{1}{2}) .
\end
{array}$$ Bounds for the term $\int_1^t\int_{\mathbb {R}^3} (|\nabla
u|^3+|\nabla u|^4)dx$ are obtained in a similar way, which in fact
is much more simple. Then applying Lemma 3.7, we can bound $\bar E$
which gives (3.31). The proof of lemma is completed.\\

Combining the results of Lemmas 3.2-3.9, we have the following
bound for $A+B_q$.\\
\\{\bf Lemma 3.10.} Assume that the hypotheses and notations of
Proposition 3.1 are in force. Then there are polynomials $\varphi_1$
and $\varphi_2$ as described in Lemma 3.9 such that for $0<t\leq
1\wedge T$,
$$A+B_q\leq
M[\varphi_1(C_0)+\varphi_2(A+B_q)],\eqno(3.32)$$ and if $T>1$ and
$1\leq t\leq T$, then
$$\bar A+B_q\leq M[\varphi_1(C_0+A(1)+B_q(1))+\varphi_2(\bar A+B_q)].\eqno(3.33)$$\\
\\{\bf Proof of Proposition 3.1.} Proposition now follows
immediately from the bounds (3.32) and (3.33) and the fact that the
functions $A,\bar A,B_q$ are continuous in time.

\section{Pointwise bounds for the density}

$\ \ \ \ $In this section we derive pointwise bounds for the density
$\rho$, bounds which are independent both of time and of initial
smoothness. This will then close the estimates of Proposition 2.1 to
give an uncontingent estimate for the functional $A$ defined in
(3.1).

We begin with two auxiliary lemmas. The first lemma is a
maximum-principle arguments applied integral curves of the velocity
field, which has been proved in Hoff [8].\\
\\{\bf Lemma 4.1} Let $(\rho,u,d)$ be as in Proposition 3.1 and
suppose that $0<c_1\leq\rho\leq c_2$ on $\mathbb R^3\times [0,T]$.
Fix $t_0\geq 0$ and define the particle trajectories
$x:[0,\infty)\times\mathbb R^3\rightarrow \mathbb R^3$ by
$$
\left\{
\begin{array}{l}
\dot x(t,y)=u(x(t,y),t),\\
x(t_0,y)=y.
\end{array}
\right.
$$
Then there is a constant $C$ depending only on $c_1$ and $c_2$ such
that if $g\in L^1(\mathbb R^3)$ is nonnegative and $t\in[0,T]$, then
each of the integrals $\int_{\mathbb R^3}g(x(t,y))dy$ and
$\int_{\mathbb R^3}g(x)dx$ is bounded by $C$ times the other.\\

Next we derive a result relating the H\"older-continuity of
$u(\cdot,t)$ to various norms appearing in the definition (2.1) of
the functional A.\\
\\{\bf Lemma 4.2.} Let $(\rho,u,d)$ be as in Proposition 3.1. Then
for $\alpha\in(0,\frac{1}{2}]$ and $t\in(0,T]$, we have
$$\begin {array}{rl}\langle u(\cdot,t)\rangle^\alpha\leq &M[\|\nabla
u(\cdot,t)\|^\frac{1-2\alpha}{2}_{L^2(\mathbb R^3)}\|\nabla
\omega(\cdot,t)\|^\frac{1+2\alpha}{2}_{L^2(\mathbb
R^3)}+(C_0+\|\nabla u(\cdot,t)\|^2_{L^2(\mathbb
R^3)})^\frac{1-2\alpha}{4}\\&\times (\|\dot
u(\cdot,t)\|^2_{L^2(\mathbb R^3)}+\|(\nabla d\cdot\triangle
d)(\cdot,t)\|^2_{L^2(\mathbb
R^3)})^{\frac{1+2\alpha}{4}}+C_0^\frac{1-\alpha}{3}].\end
{array}\eqno (4.1)$$
\\{\bf Proof.} Let $\alpha\in(0,\frac{1}{2}]$ and define $r\in(3,6]$
by $r=\frac{3}{1-\alpha}$. Then by (2.3) and (2.11), we have
$$\langle u(\cdot,t)\rangle^\alpha\leq M[\|F(\cdot,t)\|_{L^r(\mathbb R^3)}+\|\omega(\cdot,t)\|_{L^r(\mathbb R^3)}+\|(\rho-\tilde\rho)(\cdot,t)\|_{L^r(\mathbb R^3)}].\eqno(4.2)$$
By (2.1), we obtain
$$\begin {array}{rl}\|\omega(\cdot,t)\|_{L^r(\mathbb R^3)}\leq&M(\|\omega(\cdot,t)\|^\frac{6-r}{2r}_{L^2(\mathbb R^3)}
\|\nabla\omega(\cdot,t)\|^\frac{3r-6}{2r}_{L^2(\mathbb R^3)})\\
\leq &M\|\nabla u(\cdot,t)\|^\frac{1-2\alpha}{2}_{L^2(\mathbb R^3)}
\|\nabla\omega(\cdot,t)\|^\frac{1+2\alpha}{2}_{L^2(\mathbb
R^3)})\end {array}$$ and
$$\begin {array}{rl}\|F(\cdot,t)\|_{L^r(\mathbb R^3)}\leq&M(\|F(\cdot,t)\|^\frac{6-r}{2r}_{L^2(\mathbb R^3)}
\|\nabla F(\cdot,t)\|^\frac{3r-6}{2r}_{L^2(\mathbb R^3)})\\
\leq &M(\|(\rho-\tilde \rho)(\cdot,t)\|^2_{L^2(\mathbb
R^3)}+\|\nabla u(\cdot,t)\|^2_{L^2(\mathbb
R^3)})^\frac{1-2\alpha}{4}\\&\times (\|\dot
u(\cdot,t)\|^2_{L^2(\mathbb R^3)}+\|(\nabla d\cdot\triangle
d)(\cdot,t)\|^2_{L^2(\mathbb R^3)})^{\frac{1+2\alpha}{4}}]\\
\leq &M(\|C_0+\|\nabla u(\cdot,t)\|^2_{L^2(\mathbb
R^3)})^\frac{1-2\alpha}{4}\\&\times (\|\dot
u(\cdot,t)\|^2_{L^2(\mathbb R^3)}+\|(\nabla d\cdot\triangle
d)(\cdot,t)\|^2_{L^2(\mathbb R^3)})^{\frac{1+2\alpha}{4}}].\end
{array}$$ Putting the above results into (4.2) yields (4.1). Thus
the proof the lemmas is completed.\\

Now we derive the upper and lower pointwise bounds for the
density.\\
\\{\bf Proposition 4.3.} Assume that the system parameters in (1.1)
satisfy the conditions (1.5)-(1.8) and let positive numbers $N$ and
$b\leq \delta$ be given. Assume $(\rho,u,d)$ is a solution of (1.1)
on $\mathbb R^3\times[0,T]$ in the sense of Proposition 2.3 with
initial data $(\rho_0,u_0)\in H^3(\mathbb R^3)$ and $d_0\in
H_n^4(\mathbb R^3;\mathbb S^2)$ satisfying (1.10)-(1.13). Then there
are positive constants $\varepsilon,M,$ and $\theta$ depending on
the parameters and assumptions in (1.5)-(1.8), on $N$, and a
positive lower bound for $b$, such that,  if $C_0<\varepsilon$ and
$\rho(x,t)>0$ on $\mathbb R^3\times[0,T]$, then in fact
$$\underline\rho\leq\rho\leq\bar\rho\ \ on\ \ \mathbb R^3\times[0,T],\eqno (4.3)$$
and $$A(T)\leq MC_0^\theta.\eqno (4.4)$$
\\{\bf Proof.} First we choose positive numbers $\kappa$ and $\kappa'$
satisfying
$$\underline\rho<\kappa<\underline\rho+b<\bar\rho-b<\kappa'<\bar\rho.$$
Recall that $\rho_0$ takes values in
$[\underline\rho+b,\bar\rho-b]$, so that
$\rho\in[\underline\rho,\bar\rho]$ on $\mathbb R^3\times[0,\tau]$
for some positive $\tau$ by the time regularity (2.7). It then
follows from Proposition 3.1 that $A(\tau)\leq MC_0^\theta$, where
$M$ is now fixed. We shall that if $C_0$ is further restricted, then
in fact that $\kappa<\rho<\kappa'$ on all of $\mathbb
R^3\times[0,T]$, and therefore that $A(T)\leq MC_0^\theta$ as well.
We shall prove the required upper bound, the proof of the lower
bound being similar.

For $y\in\mathbb R^3$ and define the corresponding particle path
$x(t)$ by
$$
\left\{
\begin{array}{l}
\dot x(t,y)=u(x(t,y),t),\\
x(t_0,y)=y.
\end{array}
\right.
$$
Suppose that there is a time $t_1\leq \tau$ such that
$\rho(x(t_1),t_1)=\kappa'$. We may take $t_1$ minimal and then
choose $t_0<t_1$ maximal such that $\rho(x(t_0),t_0)=\bar\rho-b$.
Thus $\rho(x(t),t)\in[\bar\rho-b,\kappa']$ for $t\in[t_0,t_1]$. We
divide into two steps:\\
{\bf Step 1.} $t_0<t_1\leq T\wedge 1$

We have from the definition (1.3) of $F$ and the mass equation that
$$(\mu+\lambda)\frac{d}{dt}[\log\rho(x(t),t)-\log\tilde\rho]+P(\rho(x(t),t))-P(\tilde\rho)=-F(x(t),t).$$
Integrating from $t_0$ to $t_1$ and abbreviating $\rho(x(t),t)$ by
$\rho(t)$, etc., we then obtain
$$(\mu+\lambda)\log\rho(s)|^{t_1}_{t_0}+\int^{t_1}_{t_0}[P(s)-P(\tilde\rho) ]ds=-\int^{t_1}_{t_0}F(s)ds.\eqno (4.5)$$
We shall show that
$$\int_{t_0}^{t_1}F(s)ds\leq \tilde MC_0^\theta\eqno(4.6)$$
for a constant $\tilde M$ which depends on the same quantities as
the $M$ from Proposition 3.1 (which has been fixed). If so, then
from (3.5), we have
$$(\mu+\lambda)[\log\kappa'-\log(\bar\rho-b)]\leq -\int^{t_1}_{t_0}[P(s)-P(\tilde\rho)]ds+\tilde MC_0^\theta\leq \tilde MC_0^\theta.\eqno (4.7)$$
where the last inequality holds because $\rho(t)$ takes values in
$[\bar\rho-b,\kappa']\subset[\max\{\tilde\rho,\rho''\},\bar\rho]$,
and $P$ is increasing on $[\tilde\rho,\bar\rho]$. But (4.7) cannot
holds if $C_0$ is small depending on $\tilde M,\kappa',$ and
$\bar\rho-b$. Stipulating the smallness condition, we therefore
conclude that there is no time $t_1$ such that
$\rho(t_1)=\rho(x(t_1),t_1)=\kappa'$. Since $y\in\mathbb R^3$ was
arbitrary, it follows that $\rho<\kappa'$ on $\mathbb
R^3\times[0,\tau]$, as claimed. The proof that $\rho>\kappa$ is
similar.

To prove (4.6) we let $\Gamma$ be the fundamental solution of the
Laplace operator in $\mathbb R^3$ and apply (1.4) to write
$$\begin {array}{rl}\int^{t_1}_{t_0}F(s)ds=&\int^{t_1}_{t_0}\int_{\mathbb R^3}(\nabla_x\Gamma(x(s)-y))\rho\dot
u(y,s)dyds\\&+\int^{t_1}_{t_0}\int_{\mathbb
R^3}(\nabla_x\Gamma(x(s)-y))(\nabla d\cdot \triangle d)(y,s)dyds.
\end {array}\eqno (4.8)$$
By Lemma 4.2, the first integral on the right side of (4.8) is
bounded exactly as in Lemma 4.2 of Hoff [8]:
$$\begin {array}{rl}&\int^{t_1}_{t_0}\int_{\mathbb R^3}(\nabla_x\Gamma(x(s)-y))\rho\dot
u(y,s)dyds\\ \leq &\|\nabla\Gamma*(\rho
u)(\cdot,t_1)\|_{L^\infty(\mathbb R^3)}+\|\nabla\Gamma*(\rho
u)(\cdot,t_2)\|_{L^\infty(\mathbb R^3)}\\&+\int_0^t\int_{\mathbb
R^3}\Gamma_{x_jx_k}(x(s)-y)[u^k(x(s),s)-u^k(y,s)](\rho
u^j)(y,s)dyds\\ \leq &\tilde MC_0^\theta+\tilde
MC_0^\theta\int^1_0\langle u(\cdot,s)\rangle^\alpha ds\\ \leq
&\tilde MC_0^\theta+\tilde MC_0^\theta\int^1_0\langle
u(\cdot,s)\rangle^\alpha ds\\ \leq &\tilde MC_0^\theta+\tilde
MC_0^\theta(\int_0^1s^{-\frac{3+6\alpha}{4}}ds)^\frac{1}{2}\int^1_0(C_0+\|\nabla
u(\cdot,t)\|^2_{L^2(\mathbb R^3)}ds)^\frac{1-2\alpha}{4}\\&\times
(\int_0^1s^\frac{3}{2}(\|\dot u(\cdot,t)\|^2_{L^2(\mathbb
R^3)}+\|(\nabla d\cdot\triangle d)(\cdot,t)\|^2_{L^2(\mathbb
R^3)})ds)^{\frac{1+2\alpha}{4}}\\ \leq &\tilde MC_0^\theta,
\end {array}$$
if $\alpha<\frac{1}{6}$. Note that (3.22) holds for $q=6$, thus if
$2<r<\frac{3q}{q+3}$, by (2.4), the second integral on the right
side of (4.8) can be bounded in as
$$\begin {array}{rl}&\int^{t_1}_{t_0}\int_{\mathbb
R^3}(\nabla_x\Gamma(x(s)-y))(\nabla d\cdot \triangle d)(y,s)dyds\\
\leq &\tilde M\int_0^1\|(\nabla d\triangle d)(s)\|_{L^2(\mathbb
R^3)}+\|(\nabla d\triangle d)(s)\|_{L^r(\mathbb R^3)}ds.\\
\leq &\tilde M\int_0^1\|(|\nabla d|^4|\triangle
d|^2)(s)\|^\frac{1}{2}_{L^2(\mathbb R^3)}\||\triangle
d|^2(s)\|^\frac{1}{2}_{L^2(\mathbb R^3)}ds\\&+\int_0^1\|(\triangle
d)(s)\|_{L^3(\mathbb R^3)}\|(\nabla
d)(s)\|_{L^\frac{3r}{3-r}(\mathbb R^3)}ds\\ \leq &\tilde
MC_0^\theta.
\end {array}\eqno $$
Thus the proof of (4.6) is completed.\\
{\bf Step 2.} $1\leq t_0<t_1$.

Again by the mass equation and the definition (1.3) of $F$,
$$\frac{d}{dt}(\rho(t)-\tilde\rho)+(\mu+\lambda)^{-1}\rho(t)(P(t)-\tilde P)=(\mu+\lambda)^{-1}\rho(t)F(t).$$
Multiplying by $(\rho(t)-\tilde \rho)$ we get
$$\frac{1}{2}(\rho(t)-\tilde\rho)^2+(\mu+\lambda)^{-1}f(t)\rho(t)(\rho(t)-\tilde\rho)^2=-(\mu+\lambda)^{-1}\rho(t)(\rho(t)-\tilde\rho)F(t),\eqno(4.9)$$
where $$f(t)=(P(t)-\tilde P)(\rho(t)-\tilde\rho)^{-1}.$$ Since
$f(t)\geq 0$ on $[t_0,t_1]$, thus integrating (4.9) over
$[t_0,t_1]$, we arrive at
$$|\kappa'-\tilde\rho|^2-|\bar\rho-b-\tilde\rho|^2\leq\tilde M\int_{t_0}^{t_1}\|F(\cdot,s)\|^2_{L^\infty(\mathbb R^3)}ds.\eqno(4.10)$$
So that if we show that
$$\int_{t_0}^{t_1}\|F(\cdot,s)\|^2_{L^\infty}ds\leq \tilde MC_0^\theta. \eqno(4.11)$$
Then as in Step 1, (4.10) cannot hold if $C_0$ is sufficiently
small. Since $y\in\mathbb R^3$ was arbitrary, it follows that
$\rho<\kappa'$ on $\mathbb R^3\times[0,\tau]$, as claimed.

To prove (4.11) we apply (1.4) and (2.4) to get
$$\begin {array}{rl}\int_{t_0}^{t_1}\|F(\cdot,s)\|^2_{L^\infty(\mathbb R^3)}ds\leq &\int_{t_0}^{t_1}\|\dot u(\cdot,s)\|^2_{L^2(\mathbb R^3)}ds
+\int_{t_0}^{t_1}\|(\nabla d\cdot\triangle
d)(\cdot,s)\|^2_{L^2(\mathbb R^3)}ds\\&+\int_{t_0}^{t_1}\|\dot
u(\cdot,s)\|^2_{L^4(\mathbb R^3)}ds+\int_{t_0}^{t_1}\|(\nabla
d\cdot\triangle d)(\cdot,s)\|^2_{L^4(\mathbb R^3)}ds\\ \leq &\tilde
M C_0^\theta+\int_{t_0}^{t_1}\|\dot u(\cdot,s)\|^2_{L^4(\mathbb
R^3)}ds+\int_{t_0}^{t_1}\|(\nabla d\cdot\triangle
d)(\cdot,s)\|^2_{L^4(\mathbb R^3)}ds.\end {array}$$ The terms
integral on the right side above can be bounded as
$$\begin {array}{rl}\int_{t_0}^{t_1}\|\dot u(\cdot,s)\|^2_{L^4(\mathbb
R^3)}ds\leq &\int_{t_0}^{t_1}\|\dot
u(\cdot,s)\|^\frac{1}{2}_{L^2(\mathbb R^3)}\|\nabla\dot
u(\cdot,s)\|^\frac{3}{2}_{L^2(\mathbb R^3)}ds\\ \leq&
(\int_{t_0}^{t_1}\|\dot u(\cdot,s)\|^{2}_{L^2(\mathbb
R^3)}ds)^\frac{1}{4}(\int_{t_0}^{t_1}\|\nabla\dot
u(\cdot,s)\|^2_{L^2(\mathbb R^3)}ds)^\frac{3}{4}\\ \leq&\tilde
MC_0^\theta\end {array}$$ and
$$\begin {array}{rl}&\int_{t_0}^{t_1}\|(\nabla d\cdot\triangle d)(\cdot,s)\|^2_{L^4(\mathbb
R^3)}ds\\ \leq &\int_{t_0}^{t_1}\|(\nabla d\triangle
d)(\cdot,s)\|^\frac{1}{2}_{L^2(\mathbb R^3)}\|\nabla(\nabla
d\cdot\triangle d)(\cdot,s)\|^\frac{3}{2}_{L^2(\mathbb R^3)}ds\\
\leq& (\int_{t_0}^{t_1}\|(\nabla d\cdot\triangle
d(\cdot,s)\|^{2}_{L^2(\mathbb
R^3)}ds)^\frac{1}{4}(\int_{t_0}^{t_1}\|\nabla(\nabla d\cdot\triangle
d(\cdot,s)\|^2_{L^2(\mathbb R^3)}ds)^\frac{3}{4}\\ \leq&\tilde
MC_0^\theta,\end {array}$$ where the last inequality follows from
Proposition 3.1. Thus (4.11) is proved. The proof of Proposition is
completed.

\section{Proof of Theorem 1.2}

$\ \ \ \ $In this section, we prove Theorem 1.2 by constructing weak
solutions as limits of smooth solutions. So, we first prove the
global-in-time existence of smooth solutions with smooth initial
data which is strictly away vacuum and is only of small energy.\\
\\{\bf Proposition 5.1.} Assume that $(\rho_0,u_0,d_0)$ satisfy (2.5).
Then for any $0<T<\infty$, there exists a unique smooth solution
$(\rho,u,d)$ of (1.1)-(1.13) on $\mathbb R^3\times[0,T]$ satisfying
(2.6)-(2.9) with $T_0$ being replaced by $T$, provided the initial
energy $C_0$ satisfies the smallness condition (1.16) with
$\varepsilon>0$ being the same one as in Proposition 3.1 and
Proposition 4.3.\\
\\{\bf Proof.} The standard local existence result (Proposition
2.3) shows that the Cauchy  problem (1.1)-(1.2) admits a unique
local smooth solution $(\rho,u,d)$ on $\mathbb R^3\times[0,T_0]$. In
view of Lemma 3.2 and Proposition 4.3, we have
$$\begin {array}{rl} &A(T_0)+\sup\limits_{0\leq t\leq T_0}\int_{\mathbb R^3}(|\rho-\tilde{\rho}|^2+|u|^2+|\nabla d|^2)dx\\&\ \ \ +\int_0^{T_0}\int_{\mathbb
{R}^3}(|\nabla u|^2+|\triangle d+|\nabla d|^2d|^2)dxdt\leq MC_0,
\end {array}\eqno (5.1)$$ and
$$\underline\rho\leq\rho\leq\bar\rho\ \ on\ \ \mathbb R^3\times[0,T_0].\eqno (5.2)$$
Then the standard arguments based on the local existence results
together with the a priori bounds (5.1)-(5.2), we deduce that
$(\rho,u,d)$ is in fact the unique smooth solution of (1.1)-(1.13)
on $\mathbb R^3\times[0,T]$ for any $0<T<\infty$.\\

With the help of Proposition 5.1, we are in a position to prove
Theorem
1.2.\\
\\{\bf Proof of Theorem 1.2.} For any map $d_0\in
H_n^1(\mathbb R^3;\mathbb S^2)$, there exists $d_0^m\in
H_n^4(\mathbb R^3;\mathbb S^2)$ such that
$$\lim_{m\rightarrow \infty}\|d_0^m-d_0\|_{H^1(\mathbb R^3)}=0.$$
Let $$\rho_0^m=J_\frac{1}{m}*\rho_0,\ \ u_0^m=J_\frac{1}{m}*u_0,$$
where $J_\frac{1}{m}=J_\frac{1}{m}(x)$ is the standard mollifier.
Then $(\rho_0^m,u_0^m)\in H^3(\mathbb R^3)$ and $d_0^m-n\in
H_n^4(\mathbb R^3;\mathbb S^2)$ and the initial norm for
$(\rho_0^m,u_0^m,\nabla d_0^m)$ (i.e., the right side of (1.13) with
$(\rho_0,u_0,\nabla d_0)$ replaced by $(\rho_0^m,u_0^m,\nabla
d_0^m)$) is bounded by $C_0$. The above Proposition can be applied
to obtain a global smooth solution $(\rho^m,u^m,d^m)$ of
(1.1)-(1.13) satisfying (3.2), (4.3) and (4.4) for all $t>0$
uniformly in $m$.

In view of (2.3) and (2.11), we see from Sobolev embedding theorem
that
$$\begin {array}{rl}\langle u^m(\cdot,t)\rangle^\frac{1}{2}\leq &C\|\nabla u^m\|_{L^6(\mathbb R^3)}\\ \leq&
C(\|F^m\|_{L^6(\mathbb R^3)}+\|\omega^m\|_{L^6(\mathbb
R^3)}+\|P^m-\tilde P\|_{L^6(\mathbb R^3)}) \\ \leq& C(\tau)(1+\|\dot
u^m\|_{L^2(\mathbb R^3)}+\|\nabla d^m\cdot\triangle
d^m\|_{L^2(\mathbb R^3)})\\ \leq &C(\tau),\ \  t\geq\tau>0,\end
{array}\eqno (5.3)$$ where $F^m,\omega^m$ and $P^m$ are the
functions $F,\omega$ and $P$ with $(\rho,u,d)$ being replaced by
$(\rho^m,u^m,d^m)$.

In addition to (5.3), we also have
$$|u^m(x,t)-\frac{1}{B_{R(x)}}\int_{B_{R(x)}}u^m(y,t)dy|\leq C(\tau)R^\frac{1}{2},$$
and hence, for $0<\tau\leq t_1\leq t_2$
$$\begin {array}{rl}&|u^m(x,t_2)-u^m(x,t_1)|\\ \leq &\frac{1}{|B_{R(x)}|}|\int_{t_1}^{t_2}\int_{B_{R(x)}}|u^m(y,t)|dydt+C(\tau)R^\frac{1}{2}\\ \leq &
CR^{-\frac{3}{2}}|t_2-t_1|^\frac{1}{2}(\int_{t_1}^{t_2}\int_{B_{R(x)}}|u^m(y,t)|^2dydt)^\frac{1}{2}+C(\tau)R^\frac{1}{2}\\
\leq &
CR^{-\frac{3}{2}}|t_2-t_1|^\frac{1}{2}(\int_{t_1}^{t_2}\int_{B_{R(x)}}(|\dot
u^m(y,t)|^2+|u^m|^2|\nabla
u|^2)dydt)^\frac{1}{2}+C(\tau)R^\frac{1}{2}.\\ \leq
&C(\tau)[R^{-\frac{3}{2}}|t_2-t_1|^\frac{1}{2}+R^\frac{1}{2}].
 \end
{array}\eqno (5.4)$$ Taking $R=|t_2-t_1|^\frac{1}{4}$ in (5.4), we
get
$$|u^m(x,t_2)-u^m(x,t_1)|\leq C(\tau)|t_2-t_1|^\frac{1}{8},\ \ 0<\tau\leq t_1\leq t_2<\infty.\eqno (5.5)$$

The same estimates in (5.3) and (5.5) also hold for $d$ and $\nabla
d$. Thus, we have proved that $\{u^m\}$, $\{ d^m\}$ and $\{\nabla
d^m\}$ are uniform H\"older continuity away from $t=0$. As a result,
it follows from Ascoli-Arzela theorem that $$u^m\longrightarrow u,\
\ d^m\longrightarrow d\ \ uniformly\ on\ compact\ sets\ in\ \mathbb
R^3\times(0,\infty).\eqno (5.6)$$ Moreover, by argument in [24] (see
also [5]), we know that
$$\rho^m\longrightarrow \rho\ \ strongly\ in\ L^p(\mathbb R^3\times(0,\infty)),\ \forall p\in[2,\infty).\eqno (5.7)$$
Therefore, passing to the limit as $m\longrightarrow\infty$ by (5.6)
and (5.7) we obtain the limited $(\rho,u,d)$ which is indeed a weak
solution of (1.1)-(1.13) in the sense of Definition 1.1 and
satisfies (1.17)-(1.24).

Next we derive the large-time behavior of $(\rho,u,d)$ in (1.25).
This can be done as the ones in [4], however, for completeness we
sketch the proof here. We first deduce from the mass equation that
$$(P(\rho)-\tilde P)_t+u\cdot\nabla(P(\rho)-\tilde P)+\gamma P(\rho)div u=0.$$
Multiplying the above equation by $4(P(\rho)-\tilde P)^3$ and
integrating it over $\mathbb R^3$, we get that
$$\frac{d}{dt}\|P(\rho)-\tilde P\|^4_{L^4(\mathbb R^3)}=\int_{\mathbb R^3}(|P(\rho)-\tilde P|^4)div u-3\gamma P(\rho)(P(\rho)-\tilde P)^3div udx,$$
which, together (3.29) shows that
$$\int_1^\infty|\frac{d}{dt}\|P(\rho)-\tilde P\|^4_{L^4(\mathbb
R^3)}|dt \leq C(1+\int_1^\infty |F|^4_{L^4(\mathbb
R^3)}ds)\int_1^\infty |\nabla u|^2_{L^2(\mathbb R^3)}ds\leq C.$$ As
a result, we have
$$\|P(\rho)-\tilde P\|_{L^4(\mathbb R^3)}\longrightarrow 0\ \ as\ \ t\rightarrow \infty.$$
This, together with (3.2) and the uniform lower and upper bound of
density, shows that
$$\lim_{t\rightarrow\infty}\|\rho-\tilde\rho\|_{L^l(\mathbb R^3)}=0\eqno (5.8)$$
holds for any $l\in(2,\infty)$.

Following the argument in [4], we take a sequence
$$u^m(t,x):=u(t+m,x),$$
for all integer $m$, and $(x,t)\in\mathbb R^3\times [1,2].$ Then
from (1.24), we have
$$\lim_{m\rightarrow\infty}\int_0^1\|\nabla u^m\|_{L^2(\mathbb R^3)}=0.$$
From (1.24) again, we have
$$\|u^m\|_{H^1(\mathbb R^3)}\leq C\ uniformly\ for\ t,m.$$
Thus we arrive at
$$\lim_{m\rightarrow \infty}\|u^m\|_{L^2(\mathbb R^3)}=0\ uniformly\ for\ t.$$
That means
$$\lim_{t\rightarrow \infty}\|u(t)\|_{L^2(\mathbb R^3)}=0.\eqno (5.9)$$
For $t\geq 1$, from (2.10) and (2.11), we obtain that
$$\begin {array}{rl}\|\nabla u(t)\|_{L^6(\mathbb R^3)}\leq &C(\|F(t)\|_{L^6(\mathbb R^3)})
+\|\omega(t)\|_{L^6(\mathbb R^3)}+\|(P(\rho)-\tilde
P))(t)\|_{L^6(\mathbb R^3)}\\ \leq &C(1+\|\nabla F(t)\|_{L^2(\mathbb
R^3)} +\|\nabla\omega(t)\|_{L^2(\mathbb R^3)})\\ \leq &C(1+\| \dot
u(t)\|_{L^2(\mathbb R^3)} +\|(\nabla d\triangle d)(t)\|_{L^2(\mathbb
R^3)})\\ \leq &C.\end {array}\eqno (5.10)$$ Combining (1.24) (5.9)
and (5.10), we have
$$\lim_{t\rightarrow\infty}\|u\|_{W^{1,r}(\mathbb R^3)}=0,\eqno(5.11)$$
holds for $r\in (2,6).$

Similarly, we have
$$\lim_{t\rightarrow\infty}\|\nabla d\|_{W^{1,r}(\mathbb R^3)}=0,\eqno(5.12)$$
holds for $r\in (2,6)$. Putting (5.8) ,(5.11) and (5.12) together
gives (1.25). Thus the proof of Theorem 1.2 is
completed.\\
\\{\bf Acknowledge}\\

Guochun Wu's research was supported by China Scholarship Council
(File No. 201206310033). Zhong Tan's research was supported
Supported by National Natural Science Foundation of China-NSAF
(Grant No. 10976026) and by National Natural Science Foundation of
China-NSAF (Grant No. 11271305).

\end{document}